\documentclass[12pt]{article}

\usepackage{amsmath,amssymb,amsthm}\def\N{\mathbb{N}}\def\R{\mathbb{R}}\def\Z{\mathbb{Z}}\def\C{\mathbb{C}}\def\<{\langle}\def\>{\rangle} \def\diag{\mathrm{diag}\, }
\usepackage{xcolor}
\title{Simple cubic variance functions in $\R^n,$ Part one}

=-0.5 in
\author{Abdelhamid Hassairi\thanks{Universit\'e de Sfax, Tunisie. \texttt{abdelhamid.hassairi@fss.rnu.tu}} \ and G\'erard Letac\thanks {Institut de Math\'ematiques, Universit\'e de Toulouse, 31062 Toulouse, France. \texttt{gerard.letac@math.univ-toulouse.fr}, T\'eSA, 7 Bd de la Gare, 31500 Toulouse}}

\begin{document}
\maketitle

\begin{abstract}
The classification of natural exponential families started with the paper \cite {Morri} where Carl Morris unifies  six very familiar families by the fact that their variance functions are polynomials of degree less or equal to two. Extension of  this classification to $\R^n$ and to degree three is the subject of this paper.

\vspace{4mm}\noindent\textit{Keywords}: Actions of the group $GL(n+1,\R)$, classification of natural exponential families, multivariate Lagrange formula. variance functions. 
\end{abstract}

\section{Introduction}\subsection{Description}Typically, a natural exponential family $F$ is characterized by a convex function $k_F$ in $\R^n$
which is the logarithm of the Laplace transform of a generating measure of $F$. This function $k_F$ is called a cumulant of $F$. All the paper is based on the following idea: if we move the representative surface of $k_F$ we get locally the representative surface of a new convex function:. For instance, if $n=1$ and if we move the parabola $y=x^2$ by the $90$ degrees rotation, we get the new convex function
$y=-\sqrt{-x}.$ If this new convex function is a cumulant function $k_{F_1}$ we have created  a new exponential family $F_1$ with deep links with $F$, obtaining a powerful tool for the classification of the natural exponential families.

This paper considers the set of natural exponential families on $\R^n$, that one can identify  to the set  of their covariance matrix functions. The group $G=GL(\R, n+1)$ is operating on this set, in a natural way described in Section 3.

The group $G$ has a subgroup $G_0$, defined in \eqref{G_0},  containing the group of affinities on $\R^n$ plus an extra  parameter that one can  interpret as the Jorgensen one. For instance, in one dimension, the set $\mathcal{P}_2$ of quadratic variance functions has been  splitted by Morris \cite{Morri} into six $G_0$-orbits. However this set  $\mathcal{P}_2$ is not preserved by $G$ while it is the case of the set $\mathcal{P}_3$ of cubic  variance functions  which is splitted into four $G$-orbits.

If $F$ is a natural exponential family in $\R^n$ the covariance matrix $V_F(m)$ of the element of $F$ of mean $m$ defines the variance function $V_F$ of $F.$ The set $\mathcal{P}_2(n)$ is the set of variance functions for $\R^n$ such that all entries are polynomials in  $m=(m_1,\ldots,m_n)$ of total degree  $\leq 2.$ In other terms $V_F(m) =A(m)+B(m)+C$ where $C$ is a constant matrix, $B$ is linear and $A$ is quadratic. The description  of the elements of $\mathcal{P}_2(n)$  is an open problem, but if we  consider   the subset $\mathcal{SP}_2(n)$ of  $\mathcal{P}_2(n)$ such that $A=amm^T$ has rank one or zero, Muriel Casalis \cite{Casal}  has shown that it can be splitted into $2n+4$ $G_0$-orbits and has described them. Here again,  this set  $\mathcal{SP}_2(n)$ is not preserved by $G.$
We therefore denote by $\mathcal{SP}_3(n)$ its image by $G$. It is  a set  of  matrices $V(m)$ whose entries are polynomials of degree $\leq 3.$ The elements of $\mathcal{SP}_2(n)$ and $\mathcal{SP}_3(n)$ are respectively called the simple quadratic and simple cubic variance functions.

The  article is the first of a two parts work for the following program: to describe the $n+3$ $G$-orbits of $\mathcal{SP}_3(n)$, meeting many more or less known distributions of the literature.
This present Part one constructs only the tools of the  realization of this program: essentially the properties of the $T_g$ transformation with $g$ in $G$ in Section 3, the splitting in $n+3$ orbits in Sections 4 and 5,  the particular role of the exponential families on $\N^n$ (with $\N$ being the set on nonegative integers) in the classification with the help of Lagrange formula for several dimensions in Section 6, and a powerful extension of a result of Philippe Rouqu\`es for computing explicitely the densities of the transformations by  a certain subgroup $H$ of $G$ in Section 7.

 \subsection{Contents.}  Writing an element $g$ of $G=GL(\R, n+1)$ as a block matrix  $g=\left[\begin{array}{cc}A&b\\c^T&d\end{array}\right]$ where $d$ is a real number, $c$ and $b$ are column vectors of $\R^n$ and $A$ is a square real matrix of order $n$,  we associate  to $g$ the Moebius transform, or homography :
$h_g(m)=(Am+b)/c^Tb+d)$ sending $\R^n$ without  the hyperplane $c^Tm+d=0$ into $\R^n$, leading to an essential transformation $T_g,$
defined below on \eqref{TG}, acting on the functions from $\R^n$ to the symmetric real matrices of order $n.$
 Section 1.3 recalls the basics of natural exponential families on $\R^n$ and the fact that  the essence of such a family is both the convex function $k=\log L$ where $L$ is the Laplace transform of a generating measure of the family, and the variance function of the family. Section 2 describes the action of $G$ on the convex functions on $\R^n.$ Section 3 describes the action of $G$ via the aboved mentioned $T_g$ on an exponential family, both from the Laplace transform point of view and from variance functions, which will be our prefered choice. Section 4 is essential for the understanding of the action of $G$. It applies our new point of  view, which is  letting  $ G$ acting on the set of exponential families to the now well known case $n=1$ of the classification of the real variance functions which are polynomials with at most degree $3.$  Section 5 delineates the set of families of $\R^n$ for $n\geq 2$ which will be appropriate to a classification by the action of $G$. Section 6 is devoted to exponential families concentrated on $\N^n$ which will appear in surprizing number of times in our classification in Part two. Section 7 gives a way to compute explicitly the densities of the transformed exponential families.

\subsection{Basics for natural exponential families.} The linear space $\R^n$  is considered as a Euclidean space with its canonical scalar product. Its elements are column matrices of height $n$ and the scalar product between $c$ and $x$   can be  denoted  in two ways: $c^Tx$ or $\< c,x\>$, depending on the readability of formulas.

Given a positive measure $\mu$ on the Euclidean space $\R^n$ its Laplace transform is defined by
$$L_{\mu}(\theta)=\int_{\R^n}e^{\<\theta,x\>}\mu(dx)\leq \infty.$$
The H\"older inequality implies that the set $D(\mu)=\{\theta\in \R^n; L_{\mu}(\theta)<\infty\}$ is convex and that
$k_{\mu}=\log L_{\mu}$ defines a convex function on $D(\mu).$ In the sequel we consider the set $\mathcal{M}(\R^n)$ of measures such that the interior $\Theta(\mu)$ of $D(\mu)$  is not empty and such that $\mu$ is not concentrated on an affine hyperplane of $\R^n.$ If $\mu\in
\mathcal{M}(\R^n)$ it is standard to prove that $k_{\mu}$ is strictly convex and real-analytic on $\Theta(\mu).$ In general, a real function $k$ defined on an open convex set $\Theta\subset \R^n$ is called a cumulant  if there exists $\mu\in \mathcal{M}(\R^n)$ such that $\Theta \subset \Theta(
\mu)$ and $k=k_{\mu}$ on $\Theta.$ Let us insist on the fact that in this definition, $\mu$ is not necessarily a probability and can be unbounded.

A $\mu$ in $\mathcal{M}(\R^n)$ generates for $\theta\in \Theta(\mu)$ the following probability on $\R^n$
$$P(\theta,\mu)(dx)=e^{\<\theta,x\>-k_{\mu}(\theta)}\mu(dx)$$ while the set $F(\mu)=\{P(\theta,\mu)\ ;\theta\in \Theta(\mu)\}$ is called the natural exponential family (NEF) generated by $\mu$. Along the paper when  we will speak of exponential families, we mean NEF. The general exponential families  of the form $e^{\<a(\theta),t(w)\>-c(\theta)}\nu(dw)$ will not concern us here. Note that $F(\mu)=F(\mu_1)$ if and only if there exists $a\in \R^n$ and a real number $b$ such that $\mu_1(dx)=e^{\<a,x\>+b}\mu(dx).$ The set of such $\mu_1$ is the set of generators of $F=F(\mu)$, and in particular any element of $F$ is a generator.
From the strict convexity and real-analyticity of $k_{\mu},$ the map $\theta\mapsto m=k'_{\mu}(\theta)$ from $\Theta(\mu)$ to $\R^n$  is a diffeomorphism of $\Theta(\mu)$ onto its image $M_F\subset \R^n.$ Because $m=k'_{\mu}(\theta)=\int_{\R^n}xP(\theta,\mu)(dx),$ the set $M_F$ is called the domain of the means of $F$. If $C_F\subset \R^n$ is the closed convex support of $\mu$ then $M_F$ is a connected open set contained in the interior $C^o_F$ of $C_F$. In many cases $M_F=C^o_F$, and we say that $F$ is steep. A sufficient condition of steepness is regularity, that means that $ D(\mu)$ is open (see Barndorff- Nielsen \cite{BarNi} for a proof). In general, the problems raised by the extension of $P(\theta,\mu)$ to $\theta\in D(\mu)\setminus{\Theta(\mu)},$ although they have an important place is literature (see Chentsov \cite{Chen} or Csisz\'ar and Matus \cite{Csis}),   will not concern us either. The monograph \cite{Rio} is representative of our point of view on NEF.

The covariance matrix of $P(\theta,\mu)$ is the Hessian matrix $k''_{\mu}(\theta).$ The inverse $m\mapsto \theta= \psi_{\mu}(m)$ of the diffeomorphism $k_{\mu}'(\theta)$ maps $M_F$ onto $\Theta(\mu).$  Therefore $$V_F(m)=k''_{\mu}(\psi_{\mu}(m))$$ is the covariance matrix of the unique element of $F$ of mean $m$.  The map $m\mapsto V_F(m)$ is called the variance function of $F$. Since
$k'_{\mu}(\psi_{\mu}(m))=m,$ by derivation we get
$\psi'_{\mu}(m)k''_{\mu}(\psi_{\mu}(m))=I_n,$ which implies
$$V_F(m)=(\psi'_{\mu}(m))^{-1}.$$ This formula shows that the knowledge of $(V_F,M_F)$ gives the knowledge of $F:$ indeed $\psi'_{\mu}(m)=\psi'_{\mu_1}(m)$ imply the existence of $(a,b)\in \R^{n+1}$ such that $\mu_1(dx)=e^{\<a,x\>+b}\mu(dx).$

A last important object about a $\mu\in \mathcal{M}(\R^n)$ is its Jorgensen set $\Lambda(\mu).$ It is  defined as the set of $\lambda>0$ such that there exists $\mu_{\lambda}\in \mathcal{M}(\R^n)$ with the following property:  for all  $\theta\in \Theta(\mu)$ we have
$$L_{\mu_{\lambda}}(\theta)=(L_{\mu}(\theta))^{\lambda}.$$ A simple example is given by the Bernoulli distribution $\Lambda(\delta_0+\delta_1)=\{1,2,\ldots\}.$
In general $\Lambda(\mu)$  is an additive semigroup of $(0,\infty)$ whose structure can be quite complex (see for instance Letac, Malouche and Maurer \cite{Letama}). The particular case $\Lambda(\mu)=(0,\infty)$ is important: we say in this case that $\mu$ is infinitely divisible, and hence $P(\theta,\mu)$ is an infinitely divisible probability in the ordinary sense. For instance $\mu(dx)=1_{(0,\infty)}(x)dx$ is infinitely divisible since $\mu_{\lambda}(dx)=1_{(0,\infty)}(x)\frac{x^{\lambda-1}}{\Gamma(\lambda)}dx$ is in $\mathcal{M}(\R)$ for any $\lambda>0.$
If $F=F(\mu)$ we denote $F_{\lambda}=F(\mu_{\lambda})$ since  $\mu_1(dx)=e^{\<a,x\>+b}\mu(dx)$ implies $(L_{\mu_1}(\theta))^{\lambda}=e^{\lambda b}L_{\mu_{\lambda}}(\theta+a).$ One easily sees that $M_{F_{\lambda}}=\lambda M_F$. We have now the important formula for $\lambda \in \Lambda(\mu)$ and $m\in M_{F_{\lambda}}$
\begin{equation}\label{BENT}V_{F_{\lambda}}(m)=\lambda V_F(\frac{m}{\lambda})\end{equation} A similar formula is obtained when considering an affine isomorphism  $T$ of $\R^n$  defined by $x\mapsto A(x)+b$ where $b\in \R^n$ and $A$ is a non singular $d\times d$ matrix. Denoting by $TF$ the set of the images of the elements  of $F$ by $T$,
clearly $TF$ is a NEF with $M_{TF}=TM_F$ and we have for all $m\in M_{TF}$ the equality
\begin{equation}\label{AFFINE}
V_{TF}(m)= A^TV_F(A^{-1}(m-b))A.
\end{equation}

\section{Moving the graph of a convex function}

\subsection{The groups $G$ and $G_0$} For fixed $n$ we have defined in Section  1.1 the group $G$ as $G=GL(\R, n+1)$ of non singular real matrices of order $n+1,$ typically written as
\begin{equation}\label{G}g=\left[\begin{array}{cc}A&b\\c^T&d\end{array}\right].\end{equation}
where $A$ is a $n\times n$ matrix, $b,c\in \R^n$ and $d\in \R$ with $\det g\neq 0$. An important subgoup of $G$ is the subset $G_0$ of matrices of the form
\begin{equation}\label{G_0}g=\left[\begin{array}{cc}A&b\\0&d\end{array}\right].\end{equation} where $d>0.$ Note that $d>0$ implies that $A$ is not singular for $g\in G_0.$
\subsection{Moving the graph of a convex function on $\R^n$}
The reader could  follow the  transformation of a convex fonction below  with the help of the two examples given for $n=1$  in Section 2.3.

Consider a strictly convex and real-analytic function $k$ defined on an open convex set $\Theta\subset \R^n$ and a linear isomorphism of $\R^n\times \R$ defined by the non singular  square matrix $g$  of order $n+1$ written by blocks
$$g=\left[\begin{array}{cc}A&b\\c^T&d\end{array}\right].$$
where $A$ is a $n\times n$ matrix, $b,c\in \R^n$ and $d\in \R.$
In other terms we have for $x,x_1\in \R^n$ and $y,y_1\in \R:$

$$ \left(\begin{array}{c}x_1\\y_1\end{array}\right)=g\left(\begin{array}{c}x\\y\end{array}\right)=\left[\begin{array}{cc}A&b\\c^T&d\end{array}\right]\left(\begin{array}{c}x\\y\end{array}\right).$$

$$ \left(\begin{array}{c}x\\y\end{array}\right)=g^{-1}\left(\begin{array}{c}x_1\\y_1\end{array}\right)=\left[\begin{array}{cc}A_1&b_1\\c_1^T&d_1\end{array}\right]\left(\begin{array}{c}x_1\\y_1\end{array}\right).$$ The graph of $k$ is the set
$E=\{(x,y)\in \R^{n+1}\ ; x\in \Theta,\ y=k(x)\}.$
After transformation by $g$ it becomes the set
$$E_1=\{(x_1,y_1)\ ; c_1^Tx_1+d_1y_1=k(A_1x_1+b_1y_1)\}$$
Denote by $\Theta_1$ the interior of the image of $E_1$ by the projection  $ (x_1,y_1)\mapsto x_1.$

Suppose now that $y_1=f(x_1)$ is a real-analytic  function defined on an open subset $U$ of $\Theta_1$ such that its graph is contained in $E_1$. This last condition is saying that for $x_1\in U$ we have
$$c_1^Tx_1+d_1f(x_1)=k(A_1x_1+b_1f(x_1))$$
Taking the second differential of both sides of the above equation  implies  the equality of the two quadratic forms on $ \R^n$ which are $h\mapsto d_1 h^Tf''(x_1)h$ and
$$h^T(A_1^T+ f'(x_1)b_1^T)k''(A_1x_1+b_1f(x_1))(A_1+b_1 f'(x_1)^T)h+b_1^T k''(A_1x_1+b_1f(x_1))b_1\times h^Tf''(x_1)h$$
From the implicit function theorem, $f'$ is real-analytic and is zero  on a connected open set $V$ if and  only if $f'\equiv 0 .$ This would imply that $f$ is linear on $V$ which contradicts the fact that $k$ is strictly convex. Therefore $A_1+b_1 f'(x_1)^T$ is non zero on some  open set $V\subset U.$
Because $k$ is strictly convex $k''(A_1x_1+b_1f(x_1))$ is positive definite as well as $(A_1+b_1\otimes f'(x_1))^Tk''(A_1x_1+b_1f(x_1))(A_1+b_1 f'(x_1)^T).$
As a consequence if $$s(x_1)=d_1-k''(A_1x_1+b_1f(x_1))b_1^T $$ then the matrix $s(x_1)f''(x_1)$ is positive definite on $V$.
Therefore if $U_1\subset V$ is a convex open set such that $$s(x_1)>0$$ for $x_1\in U_1$, we have $f$ strictly convex on $U_1.$

This condition $s(x_1)>0$ for all $x_1\in U_1$ is quite important and is in fact equivalent to assume that $f$ is strictly convex on $U$. Actually if $g$ is the simple symmetry
$$g=\left[\begin{array}{cc}I_n&0\\0&-1\end{array}\right]$$
$f=-k$ is concave and of no interest for our purposes. Indeed we intend to use this construction by $g$ of a new convex function $f(x_1)$ from an old one $k(x)$ when $k(\theta)=k_{\mu}(\theta)$ is a cumulant function in order to obtain sometimes a new cumulant function $f(\theta_1)=k_{\mu_1}(\theta_1).$ We say 'sometimes' since this is not  always the case: for instance, consider for $n=1$ the cumulant function of the  ordinary Bernoulli distribution $k(\theta)=\log{\frac{1+e^{\theta}}{2}}$ and the symmetry
$$g=\left[\begin{array}{cc}0&1\\1&0\end{array}\right]$$ leading to $k_1(\theta_1)=\log(2e^{\theta_1}-1).$ Clearly $k_1$ while convex is not a cumulant function since $e^{k_1}$ is the Laplace transform of the non positive measure $2\delta_1-\delta_0.$

\subsection{Examples for $n=1$}

\vspace{4mm}\noindent\textbf{Example 2.1.} \textsc{From Gaussian laws to inverse Gaussian laws:}
Let us illustrate the above construction in the particular case $n=1$  with the convex function $k_{N(0,1)}(\theta)=\theta^2/2$. Its graph is a parabola. Consider  for real $r$ the linear transformation with matrix :
\begin{equation}\label{ROTATION}g_r=\left[\begin{array}{rr}\cos r&-\sin r\\\sin r&\cos r \end{array}\right]\end{equation}
We obtain, where we have replaced $x_1$ by $\theta_1$ for getting a familiar result:
\begin{equation}\label{PARABOLA}y_1=k_{\mu_1}(\theta_1)=\frac{1}{\sin^2 r}(\cos r-\theta_1\cos r \sin r \, -\sqrt{\cos^2 r-2\theta_1\sin r\,  })\end{equation} where $\theta_1\sin r<\cos^2 r.$ Its graph is a part of a parabola. The particular case $r=\pi/2$ gives in particular $k_{\mu_1}(\theta_1)=-\sqrt{-2\theta_1}$ where $\theta_1<0$. Since \begin{equation}\label{STABLE}\int_0^{\infty}e^{\theta_1 x-\frac{1}{2x}}x^{-3/2}\frac{dx}{\sqrt{2\pi}}=e^{-\sqrt{-2\theta_1}}\end{equation} we see that $\mu_1$ exists for $r=\pi/2$. Actually it exists for all values of $r$ since \eqref{PARABOLA} leads to a Laplace transform of a measure which is the image by an affine transform of an inverse Gaussian measure  proportional to $e^{-ax-\frac{b}{x}}x^{-3/2}1_{(0,\infty)}(x)dx.$

\vspace{4mm}\noindent\textbf{Example 2.2.} \textsc{Circle law:} Consider the probability $\mu$ such that
$$k_{\mu}(\theta)=(1-\sqrt{1-\theta^2}),\ \ \ \Theta(\mu)=(-1,1).$$
The existence of $\mu$ is granted by the formula
\begin{equation}\label{ESM}\int_0^{\infty}e^{\frac{\theta^2}{2}y}\nu(dy)=e^{(1-\sqrt{1-\theta^2})}\end{equation}
where $$\nu(dy)=\frac{1}{\sqrt{2\pi}}y^{-3/2}e^{1-\frac{y}{2}-\frac{1}{2y}}1_{(0,\infty)}(y)dy$$ is an inverse Gaussian distribution. Clearly $\mu$ is infinitely divisible. The graph of $k_{\mu}$ is a half of  a circle. Apply  now the transformation \eqref{ROTATION} for $0<r<\pi/2.$ The new convex function $k_{\mu_1}(\theta_1)$ has a graph which is an arc of the same circle. An important point is that arc is strictly contained in $\Theta(\mu_1).$ Because of the real-analyticity $\mu=\mu_1.$


\section{The actions of $G$ and $G_0$ on the exponential families}
In this section we explain how an element of $G_0$ and $G$ may transform
an exponential family $F$ into another one. This action can be seen in two ways: either by acting on the cumulant function of a generating measure of $F$, or by acting directly on the variance function of $F.$ The word 'action' should not been quite taken in the sense of group theory, since it is an action on functions (cumulant or variance), but
as we have seen in Section 2, this transformation does not always yield the cumulant of a positive measure, or similarly a function which is a variance function. The action on the variance is simpler to visualize, but ignores  the hard work of the interpretation of  the new variance by its Laplace transform and by a generating measure. The action on variance functions is described in Sections 3.2 and 3.3, the action on the cumulant functions is described in Section 3.4.

\subsection{The homography associated to an element of $G$}
Given $g=\left[\begin{array}{cc}A&b\\c^T&d\end{array}\right]$ in $G$ we consider the map $x\mapsto  h_g(x)$ from $\R^n$ into itself (called homography) defined by

\begin{equation}\label{HOMOGRAPHY}
h_g(x)=\frac{Ax+b}{c^Tx+d}\end{equation}
Since this does not make sense for $x$ in the hyperplane $H_g=\{x; \ c^Tx+d =0\} $ we mention that a proper way to deal with this transformation would be to inbed $\R^n$ in its natural projective space. However it will not be really useful for our purposes: we will simply avoid to use $h_g$ on the hyperplane $H_g.$ A simple calculation shows that
\begin{equation}\label{COMPOSITION}h_{g_1}\circ h_{g}=h_{g_1g}.\end{equation}

\vspace{4mm}\noindent\textbf{Proposition 3.1.} For $x\not \in H_g$, the differential $h'_g(x)$ is the endomorphism of $\R^n$ equal to
$$\Delta\mapsto h'_g(x)(\Delta)=
\frac{1}{(c^Tx+d)^2}[(c^Tx+d)A \Delta-(Ax+b)c^T \Delta]$$ and it is non singular. Written as a square matrix of order $n$ we have $$h'_g(x)=
\frac{1}{(c^Tx+d)^2}[(c^Tx+d)A -(Ax+b)c^T ]$$
In particular \begin{equation}\label{DERIVE}h'_{g_1g}(x)=h'_{g_1}(h_g(m))\times h'_g(m).\end{equation}
\vspace{4mm}\noindent\textbf{Proof.} The calculation of $h'_g(x)$  is done by considering  $h_g(x+\Delta)-h_g(x)$ and by using the definition of the differential. For seeing that  $h'_g(x)$ is non singular, assume that there exists $\Delta\in \R^n\setminus \{0\}$ such that $h'_g(x)(\Delta)=0.$ One easily sees that
\begin{equation}\label{MIRACLE}A\Delta=h_g(x)c^T\Delta.\end{equation} Introduce for such a $\Delta$ the vector $u$ of $\R^{d+1}$ defined by
$$u=\left(\begin{array}{c}(c^Tx+d)\Delta-xc^T\Delta\\-c^T\Delta\end{array}\right)$$ Observe that $u=0$ would imply $\Delta=0$ since the number $c^Tx+d$ is not zero. However, using \eqref{MIRACLE}

$$g(u)=\left(\begin{array}{c}(c^Tx+d)A\Delta-Axc^T\Delta -bc^T\Delta\\c^T(c^Tx+c)\Delta -c^Txc^T\Delta\-dc^T\Delta\end{array}\right)=\left(\begin{array}{c}(c^Tx+d)h_g(x)c^T\Delta-Axc^T\Delta -bc^T\Delta\\0\end{array}\right)=0$$ Since $u\neq 0$ and $g(u)=0$ this contradicts the fact that $g$ is not singular and proves the proposition. Formula \eqref{DERIVE} is clear.

\subsection {The action of $G$ on variance functions}
Here is the  fundamental definition  of $T_g:$

\vspace{4mm}\noindent\textbf{Definition 3.1}. Let $g$ as defined by \eqref{G}, let $U$ be  a non void open subset of $\R^n$ such that \begin{equation}\label{BASE}c_Tm+d>0\end{equation} for all $m\in U$ and let $m\mapsto V(m)$ be a map from $U$ to the cone on positive definite matrices. We denote

    \begin{equation}\label{TG}\boxed{T_g(V)(m)=\frac{1}{c^Tm+d}\times (h'_g(m))^{-1}V(h_g(m)) (h'_g(m)^T)^{-1} }\end{equation}

Clearly the matrix $T_g(V)(m)$ is positive definite. A germane  example (see Letac and Mora \cite{Letmo}) is $$n=1,\ U=(0,\infty),\ g=\left[\begin{array}{cc}0&1\\1&0\end{array}\right], \ \ h_g(m)=\frac{1}{m},\ T_g(m)=m^3V(\frac{1}{m}).$$

\vspace{4mm}\noindent\textbf{Proposition 3.2.} $T_{g}\circ T_{g_1}=T_{g_1 g}.$

\vspace{4mm}\noindent\textbf{Proof.} The proof  is a tedious verification:

 $$g_1g=\left[\begin{array}{cc}A_1&b_1\\c_1^T&d_1\end{array}\right]\left[\begin{array}{cc}A&b\\c^T&d\end{array}\right]= \left[\begin{array}{cc}A_1A+b_1c^T&A_1b+b_1d\\c_1^TA+d_1c^T&c_1^Tb+d_1d\end{array}\right]$$

\begin{eqnarray}\nonumber&&T_g(T_{g_1}(V)(m)=\frac{1}{c^Tm+d}(h_g'(m))^{-1}T_{g_{1}}(V)(h_g(m))(h'_g(m)^T)^{-1}\\
\nonumber&=&\frac{1}{c^Tm+d}\times \frac{1}{c_1^Th_g(m)+d_1}(h_g'(m))^{-1}(h_{g_1}'(h_g(m)))^{-1}V(h_{g_1}\circ h_g(m))(h_{g_1}'(h_g(m))^T)^{-1}(h'_g(m)^T)^{-1}\\\label{COMPO1}&=&\frac{1}{c^Tm+d}\times \frac{1}{c_1^Th_g(m)+d_1}\times (h_{g_1g}'(m)^{-1}V(h_{g_1g}(m)))((h_{g_1g}'(m)^T)^{-1}\\\label{COMPO2}&=&\frac{1}{(c_1^TA+d_1c^T)m+c_1^Tb+d_1d}\times (h_{g_1g}'(m)^{-1}V(h_{g_1g}(m)))((h_{g_1g}'(m)^T)^{-1}=T_{g_1 g}(V)(m).
\end{eqnarray}
Line \eqref{COMPO1} is the consequence of \eqref{DERIVE} and \eqref{COMPO2} results from  a little calculation.

\subsection {The action of $G_0$ on variance functions. The $G_0$-orbits. }
Let us examine $T_g$  for  two typical elements $g$ of $G_0.$

\begin{itemize}

\item $g=\left[\begin{array}{cc}A&b\\0&1\end{array}\right]$,  $h_g(m)=Am+b$, $ h'_g(m)=A$, $T_g(V)(m)=(A^T)^{-1}V(Am+b)A^{-1}$ while $T_{g}(V)$ gives formula \eqref{AFFINE}. We will call such a $g$ an \textit{affine element} of $G.$

\item $g=J_{\lambda}=\left[\begin{array}{cc} I_n&0\\0&\lambda\end{array}\right]$, $h_g(m)=m/\lambda$,  $ h'_g(m)=I_n/\lambda$ and
$T_g(V)(m)=\lambda V(m/\lambda)$
which is formula \eqref{BENT}. We will call  $J_{\lambda}$  a \textit{Jorgensen element} of $G.$
\end{itemize}

This leads to two  possible  decompositions into a product of affine and Jorgensen elements for any  $g\in G_0$:
\begin{equation}\label{DECGZERO}
\left[\begin{array}{cc}A&b\\0&d\end{array}\right]=\left[\begin{array}{cc} A&b/d\\0&1\end{array}\right]\left[\begin{array}{cc} I_n&0\\0&d\end{array}\right]=\left[\begin{array}{cc}I_n&0\\0&d\end{array}\right]\left[\begin{array}{cc} A& b\\0&1\end{array}\right].
\end{equation}

\vspace{4mm}\noindent\textbf{Comment.} Let us clarify a delicate point. We see that the knowledge of $h_g$ does not give a complete  knowledge of $g$ since $h_{\alpha\,  g}=h_g$  where $\alpha\neq 0$ is a scalar. Therefore, we really need a complete knowledge of $g$  with its $(n+1)^2 $ parameters to define $T_g$. For example, taking  $n=1$ the transformation $V(m)\to a^2V(m/a)$
gives the variance function of the image of the exponential family by $x\mapsto ax.$
The transformation $V(m)\mapsto \lambda V(m/\lambda)$  is the Jorgensen transformation and is of a quite different nature.
The action is from $g=\left[\begin{array}{cc}a&b\\0&d\end{array}\right]$ is not from the affinity $h_g(m)=\frac{am+b}{d}$ since
$T_g(V)(m)=\frac{d}{a^2}V\left((\frac{am+b}{d}\right):$ we have to resist to the tentation of taking $d=1$ for simplification.

\vspace{4mm}\noindent\textbf{Definition 3.2}. Let $V$  the variance function of some exponential family $F$ in $\R^n.$ The $G_0$-orbit is the set of all variance functions $V_1$ such that $T_g(V)(m)=V_1(m)$ for some $g\in G_0$ and some non void open subset of the domain of the means $M_1$ for $V_1$.

\vspace{4mm}\noindent
Real-analyticity of the variance functions implies the uniqueness of $V_1.$ Since there is  a one to one correspondence between the set of variance functions and the set of exponential families we  speak also of the  $G_0$-orbit of an exponential family.
For instance Carl Morris \cite{Morri} presents six groups of exponential families corresponding to the six $G_0$-orbits of the variances $1, m, m^2, m-m^2, m+m^2, m^2+1$ which are essentially the $G_0$-orbits
of normal, Poisson, Gamma, binomial, negative binomial and hyperbolic. We say 'essentially', since  the elements of a $G_0$-orbit include the affinities while an affinity of a Poisson distribution for instance has a small practical interest. We will comment later about a generalisation to $\R^n$ of the Morris classification made by Muriel Casalis \cite{Casal} which is the source of the present paper.

Here is another example of $G_0$-orbit.

\vspace{4mm}\noindent\textbf{Example 2.2, continued.} $k_{\mu}(\theta)=1-\sqrt{1-\theta^2}$ for $|\theta|<1$ implies that the corresponding variance function is
$V(m)=(m^2+1)^{3/2}.$  The action of $g=\left[\begin{array}{cc}A&b\\0&d\end{array}\right]$ on $V$ gives
$$T_g(V)(m)=((am+b)^2+d^2)^{3/2}\frac{d^4}{a^2}.$$
Changing a bit our notations for simplicity,
we denote by $C$ the cone of $(a,b,c)\in \R^3$ such that $a,c,ac-b^2$ are positive. In other terms $(a,b,c)$ belongs to $C$
if and only if $am^2+2bm+c>0$ for all $m.$ Observe that in this case $$V_{F(a,b,c)}(m)=(am^2+2bm+c)^{3/2}$$ is the variance function of some exponential family $F(a,b,c)$ such that $M_F=\R.$ Therefore the NEF $F(a,b,c)$ of variance function $(am^2+2bm+c)^{3/2}$ is the image of $V$ by the affine map
$$y=h(x)=\frac{\sqrt{ac-b^2}}{a}x-\frac{b}{a}$$ completed by a suitable Jorgensen transformation, which is licit since $\mu$ is infinitely divisible.
 The set $ \{F(a,b.c); (a,b,c)\in C\}$ is  a $G_0$-orbit.

\subsection{The action of $G$ on the Laplace transforms}

\vspace{4mm}\noindent\textbf{Definition 3.3.} If $U$ is an open subset of $\R^n$ and if $g=\left[\begin{array}{cc}A&b\\c^T&d\end{array}\right]\in G$ we denote
\begin{equation}\label{UG}U_g=\{m\in \R^n; \ c^Tm+d>0,\ h_{g}(m)\in U\}.\end{equation}
In this section, we link Section 2  to the action of $G$ on variance functions defined by \eqref{TG} with the following theorem. (Note that here $\lambda$ denotes a vector of $\R^n$ and is not a Jorgensen parameter).

\vspace{4mm}\noindent\textbf{Theorem 3.3.} Let $F=F(\mu)$ and $F_1=F(\nu)$ be two exponential families on $\R^n.$ Then  there exists $g=\left[\begin{array}{cc}A&b\\c^T&d\end{array}\right]$ in $G$ such that $F=T_g(F_1)$  if and   only if there exists an open subset $U$ of $\Theta(\mu)$, a map $\lambda \mapsto \theta(\lambda)$ from $U$ to $\Theta(\nu)$ and a point $\left[\begin{array}{c}\lambda_0\\-k_0\end{array}\right]\in \R^{n+1}$ such that
\begin{equation}\label{RECIPROCITY}
\left[\begin{array}{c}\lambda\\-k_{\mu}(\lambda)\end{array}\right]=\left[\begin{array}{cc}A^T&c\\b^T&d\end{array}\right]\left[\begin{array}{c}\theta(\lambda)\\-k_{\nu}(\theta(\lambda))\end{array}\right]
+\left[\begin{array}{c}\lambda_0\\-k_0\end{array}\right]
\end {equation}

\vspace{4mm}\noindent\textbf{Proof.} Suppose that $F=T_g(F_1).$
Consider the open set $M=M_F\cap (M_F)_g.$ It is not empty from the hypothesis and 	definition of $T_g.$ Define $U=\psi_{\mu}(M)$ where $\psi_{\mu}(m)=\theta$ is the reciprocal of $\theta\mapsto m=k'_{\mu}(\theta)$ as defined in Section 1. The implicit function theorem implies that $\psi_{\mu}$ is real-analytic and is an open application. Therefore $U$ is open.  It is contained in $ \Theta(\mu)$ and we define the map from $U$ to $\Theta(\nu)$ as
$$\theta : \lambda\mapsto \theta(\lambda)=\psi_{\nu}(h_g(k'_{\mu}(\lambda))).$$
As a consequence \begin{equation}\label{THETALAMBDA1} k'_{\nu}(\theta(\lambda))=h_g(k_\mu'(\lambda))\end{equation}
Taking derivative:

\begin{equation}\label{THETALAMBDA2} k''_{\nu}(\theta(\lambda))\theta'(\lambda)=h'_g(k_\mu'(\lambda))k''_{\mu}(\lambda)\end{equation}
The fact that $V_{F_{\mu}}(m)=T_g(V_{F_{\nu}})(m)$  for all $m\in M$ is saying that we have
\begin{equation}\label{THETALAMBDA3}k''_{\mu}(\lambda)=\frac{1}{c^Tk'_{\mu}(\lambda)+d}\times (h'_g(k'_{\mu}(\lambda)))^{-1}V_{F_{\nu}}(h_g(k'_{\mu}(\lambda))) (h'_g(k'_{\mu}(\lambda))^T)^{-1}\end{equation} by replacing $m$ by $m=k'_{\mu}(\lambda).$
Comparing \eqref{THETALAMBDA2} and \eqref{THETALAMBDA3}  and using \eqref{THETALAMBDA1} we get

$$k''_{\nu}(\theta(\lambda))\theta'(\lambda)=\frac{1}{c^Tk'_{\mu}(\lambda)+d}\times k''_{\nu}(\theta(\lambda)) (h'_g(k'_{\mu}(\lambda))^T)^{-1}.$$ which implies
$$\theta'(\lambda)=\frac{1}{c^Tk'_{\mu}(\lambda)+d}\times (h'_g(k'_{\mu}(\lambda))^T)^{-1}, $$ rewritten as
\begin{equation}\label{THETALAMBDA4}
 I_n=(c^Tk'_{\mu}(\lambda)+d)\times  (h'_g(k'_{\mu}(\lambda))^T)\theta'(\lambda)\end{equation}
Using Proposition 3.1 and then \eqref{THETALAMBDA1} equation \eqref{THETALAMBDA4} becomes
\begin{eqnarray}\nonumber I_n&=&\left[A^T-c\frac{(Ak'_{\mu}(\lambda)+b)^T}{c^Tk'_{\mu}(\lambda)+d}\right]\theta'(\lambda)=\left[A^T-ch_g(k'_{\mu}{\lambda})^T\right]\theta'(\lambda)\\\label{THETALAMBDA5}&=&\left[A^T-c k'_{\nu}(\theta(\lambda))^T\right]\theta'(\lambda)\end{eqnarray} Integrating \eqref{THETALAMBDA5} there exists $\lambda_0$ such that
$$\lambda-\lambda_0=A\theta(\lambda)-ck_{\nu}(\theta(\lambda)).$$
Similarly for revealing $k_0$ we rewrite \eqref{THETALAMBDA1} as
$$(A-k'_{\nu}(\theta(\lambda))c^T) k'_{\mu}(\lambda)=k'_{\nu}(\theta(\lambda))-b$$
Transposing  the last line, multiplying by $\theta'(\lambda)$, applying \eqref{THETALAMBDA5} and integrating, we obtain as desired  the existence of a real number $k_0$ such that $$k_{\mu}(\lambda)-k_0=b\theta(\lambda) -k_{\nu}(\theta(\lambda)).$$ The only if part of the proposition is proved.

\vspace{4mm}\noindent Conversely, suppose that  \eqref{RECIPROCITY} holds. For simplification, we denote

$$m=k'_{\mu}(\lambda), \ M=k'_{\nu}(\theta(\lambda)).$$

Taking the first and the second differential in \eqref{RECIPROCITY} we get, using $\Delta$ and $\Delta_1$ in $\R^n$
\begin{eqnarray}\label {D1}\Delta &=&(A^T -cM^T)\theta'(\lambda)\Delta ,
\\ \label {D2}
-m^T\Delta&=&(b -dM^T)\theta'(\lambda)\Delta
\\\label {D3}
0&=&(A^T-cM^T)\theta''(\lambda)(\Delta,\Delta_1)-c k''_{\nu}(\theta(\lambda))(\theta'(\lambda)\Delta,  \theta'(\lambda)\Delta_1 )
\\\label {D4}
-k''_{\nu}(\lambda)(\Delta,\Delta_1)&=&(b-dM)^T\theta''(\lambda)(\Delta,\Delta_1)-d k''_{\nu}(\theta(\lambda))(\theta'(\lambda)\Delta,  \theta'(\lambda)\Delta_1 )
\end{eqnarray}

Since \eqref{D1} says that $\theta'(\lambda)=(A^T -cM^T)^{-1}$, carrying in \eqref{D3} gives

$$\theta''(\Delta,\Delta_1)=\theta'(\lambda)ck''_{\nu}(\theta(\lambda))(\theta'(\lambda)\Delta,  \theta'(\lambda)\Delta_1 )$$
We bring this in \eqref{D4}:
$$-k''_{\nu}(\lambda)(\Delta,\Delta_1)=(b-dM)^T\theta'(\lambda)ck''_{\nu}(\theta(\lambda))(\theta'(\lambda)\Delta,  \theta'(\lambda)\Delta_1 )-d k''_{\nu}(\theta(\lambda))(\theta'(\lambda)\Delta,  \theta'(\lambda)\Delta_1 )$$
We use \eqref{D2} for simplification :

\begin{eqnarray}\nonumber k''_{\nu}(\lambda)(\Delta,\Delta_1)&=&m^T\, c \, k''_{\nu}(\theta(\lambda))(\theta'(\lambda)\Delta,  \theta'(\lambda)\Delta_1 )+d k''_{\nu}(\theta(\lambda))(\theta'(\lambda)\Delta,  \theta'(\lambda)\Delta_1 )\\ \label{D5}&=&
(m^T c +d) \, k''_{\nu}(\theta(\lambda))(\theta'(\lambda)\Delta,  \theta'(\lambda)\Delta_1 )
\end{eqnarray}
which can be rewritten as
$$\Delta^Tk''_{\nu}(\lambda)\Delta_1=(m^T c +d)\Delta^T\theta'(\lambda)^Tk''_{\nu}(\theta(\lambda))\theta'(\lambda)\Delta_1$$
$$=(m^T c +d)\theta'(\lambda)^TV_{F_{\mu}}(M)\theta'(\lambda)$$

Formulas \eqref{D1} and \eqref{D2} imply $$ k''_{\nu}(\theta(\lambda))=h_g(m)$$ and
\begin{eqnarray}\nonumber\theta'(\lambda)^T&=&(A-Mc^T)^{-1}=(A-h_d(m)c^T)^{-1}\\\label{D6}&=&\left(\frac{A(c^Tm+d)-(Am+b)c^T}{c^Tm+d}\right)^{-1}=\frac{1}{c^Tm+d}\times (h'_g(m))^{-1}.\end{eqnarray}
we get as desired $$V_{F_{\nu}}(m)=T_g(V_{F_{\mu}})(h_g(m)).$$  Theorem 3.3 is proved.

Let us now give  a \textit{necessary} algebraic condition on functions $m\mapsto V(m)$ to be a variance function.

\vspace{4mm}\noindent\textbf{Proposition 3.4.} Let $M$ be an open set of $\R^n$ and let $S$ be the set of continuously differentiable functions $V$ from $M$ to the set  of symmetric  matrices of order $n$ such that the following bilinear map from $\R^n\times \R^n$ to $\R^n$  defined by
$$(\Delta,\Delta_1)\mapsto f(\Delta,\Delta_1)=V'(m)(V(m)\Delta, \Delta_1)$$
satisfies $f(\Delta,\Delta_1)= f(\Delta_1,\Delta).$ Then \begin{enumerate}
\item Any variance function defined on $M$ is in $S$
\item $S$ is stable by $T_g$ for any $g$ such that $M_g$ is not empty.
\end{enumerate}

\vspace{4mm}\noindent\textbf{Proof.} 1) Let $\mu$ be a generating measure of  the exponential family $F$. As mentioned in Section 1 we have
$$\psi'_{\mu}(m)=V_F(m)^{-1}.$$ Since $\psi_{\mu}(m)\in \R^n$ and since $\psi''_{\mu}$ is Hessian we can consider $\psi''_{\mu}$ as a symmetric bilinear function of two variables valued in $\R^n.$
Taking derivative with respect to $m$ of $\psi'_{\mu}(m)x=V_F(m)^{-1}x\in \R^n$ where $x\in \R^n$ we get
$$ \psi''_{\mu}(m)(x,y)=-V_F(m)^{-1}(V'_F(m)y)V_F(m)^{-1} x\in \R^n.$$

The above line in symmetric in $x$ and $y$. Hence
$$(V'_F(m)x)V_F(m)^{-1} y=(V'_F(m)y)V_F(m)^{-1} x\in \R^n,$$ which can be rewritten as
$$(V'_F(m)(x,V_F(m)^{-1} y)=(V'_F(m)(y,V_F(m)^{-1} x)\in \R^n,$$
and we have just to coin $\Delta=V_F(m)^{-1}x$ and $\Delta_1=V_F(m)^{-1}y$ to get the result.

\vspace{4mm}\noindent 2)
 From \eqref{TG} we compute the differential of$$ m\mapsto T_g(V)(m)=\frac{1}{c^Tm+d}\times (h'_g(m))^{-1}V(h_g(m)) (h'_g(m)^T)^{-1}$$
To perform this computation we use the auxiliary result
\begin{equation}\label{DDH}  ((h_{g} ' {(m)})^{-1} )'(t)  = -(h_{g} '(m) )^{-1}h_g''(m) (t)(h_{g} ' (m))^{-1} \end{equation}
and we get
\begin{eqnarray}
\label {H1}T_g(V)'(m)(t)&=&-\frac{c^T t}{(c^Tm+d)^2}\times (h'_g(m))^{-1}V(h_g(m)) (h'_g(m)^T)^{-1}\\
\label {H2}&&-(h_{g} '(m) )^{-1}h_g''(m) (t)(h_{g} ' (m))^{-1}V(h_g(m)) (h'_g(m)^T)^{-1}\\
\label {H3}&&+(h_{g} '(m) )V'(h_{g} (m) )((h_{g} '(m) )(t)(h'_g(m)^T)^{-1}\\
\label {H4}&&-(h_{g} '(m) )^{-1}V'(h_{g} (m) )((h_{g} '(m))^T )^{-1}h_g''(m) (t)((h_{g} ' (m))^T)^{-1}\end{eqnarray}

Now we replace $t$ by $T_g(m)\Delta$ ,  and we verify that each line \eqref{H1}, \eqref{H2}, \eqref{H3}, \eqref{H4}
will fit with the desired symmetry between $\Delta$ and $\Delta_1$  in the expression
$$[(T_g V)'(m)][(T_g V)(m)(\Delta)]\Delta_1.$$
For lines \eqref{H1} and \eqref{H4} symmetry comes from the fact the $V(m)$ is a symmetric matrix. For line \eqref{H2} symmetry comes from the fact that $h_g''(m)(\Delta, \Delta_1)=h_g''(m)(\Delta_1, \Delta)$ and  for line \eqref{H3} symmetry comes from the fact that $V$ is in $S.$
\section{The  $G$-orbits of the   Morris Mora class for $\R$}
\subsection {The Morris class for $\R$ and its six $G_0$-orbits}
We have already mentioned  the Morris class $\mathcal{M}_2(\R)$ of exponential families as the set of exponential families such that their variance function is a polynomial of degree $\leq 2$. In particular, it contains the following six variance functions $(V_F, M_F) $
\begin{eqnarray*}(1,\R),&&Normal \\ (m,(0,\infty),&& Poisson\\ (m^2,(0,\infty),&& Gamma \\ (m-m^2,(0,1)),&& Binomial\\  (m+m^2,(0,\infty)),&& Negative Binomial\\  ( m^2+1,\R),&& Hyperbolic\end{eqnarray*}
An important point is the fact that  each of the six $G_0$-orbits   that they generate are disjoint since the action of $G_0$ on a variance function respects 1) the degree, 2) the number of real zeros, complex zeros and their multiplicities 3) the sign of the coefficient of $m^2$. Criteria 1) isolate Normal and Poisson, criteria 2) separates Gamma, Hyperbolic  and the pair Binomial and Negative Binomial, Critera 3) separates Binomial and Negative binomial.
Note that for  $g=\left[\begin{array}{cc}a&b\\0&d\end{array}\right]$ we have $d>0$ from  \eqref{BASE}, but $a$ can be negative. For instance if $a=-1,b=0,d=1$ applying  $T_g$ to Gamma send it to $(m^2,(-\infty,0)$ which also belongs to the $G_0$-orbit of Gamma.

\subsection {The Morris Mora class in $\R$ and its four $G$-orbits.} Let us observe that the Morris class is unstable by $G.$ Taking $g=\left[\begin{array}{cc}0&-1\\1&0\end{array}\right]$ applied through $T_g$ to the constant variance function $V(m)=1$ gives $T_g(V)(m)=m^3$ with $(0,\infty) $ as domain of the means, getting a polynomial of degree three. For $g=\left[\begin{array}{cc}a&b\\c&d\end{array}\right]$ we get if $V$ is a polynomial of degree $\leq 3$
$$T_g(V)(m)=(cm+d)^3V\left(\frac{am+b}{cm+d}\right)$$ which is also a polynomial of degree $\leq 3.$ Therefore in order to understand the action of $G$ on the Morris class, we have to know all the variance functions which are polynomials of degree $\leq 3$. This description has been done by Marianne Mora in her PhD Thesis, and we call this  whole set of variance the Morris Mora class $\mathcal{M}_3(\R)$. We are not going to describe $\mathcal{M}_3(\R)$ here, since this is detailed in the paper Letac and Mora \cite{Letmo}.

Splitting $\mathcal{M}_3(\R)$ in $G_0$-orbits is not interesting since the class has four parameters and $G_0$ has only three. In the other hand a splitting of $\mathcal{M}_3(\R)$ in $G$-orbits is reasonable since $G$ has four parameters. Actually an  action of $G=GL(\R, 2)$ on the set of polynomials of degree $\leq n$ by
$$P\mapsto (cm+d)^nP\left(\frac{am+b}{cm+d}\right)$$  is of interest to geometers. The case $n=3$ is mentioned in Arnold et \textit {al} \cite{Arvag}  page 157 where the  set of polynomials of degree $\leq 3$ is splitted  in four classes generated by the following polynomials
$$ X^3, \ X^2, \ X(X+1), \ X^2+1$$ and each element of $\mathcal{M}_3(\R)$ belongs to one of these four orbits: for instance normal class belongs to $X^3$, Poisson and Gamma belong to $X^2$,  Hyperbolic belongs to $X^2+1$ and Binomial and Negative Binomial belong to  $X(X+1).$

As a result, the Morris Mora class has four $G$-orbits.

\section {The  Casalis class $\mathcal{SM}_2(\R^n)$, the Morris Mora class  $\mathcal{SM}_3(\R^n)$}
The most natural generalisation of the Morris class to $\R^n$ is the set $\mathcal{M}_2(\R^n)$, of variance matrices of the form
$$V(m)= A(m)+B(m)+C$$where $A(m)$ is a symmetric matrix of order $n$ whose entries are homogeneous quadratic polynomials in $m=(m_1,\ldots,m_n)$, where $B(m)=B_1m_1+\cdots+B_nm_n$ and  $B_1,\ldots,B_n,C$ are  symmetric matrices of order $n.$ A classification of $\mathcal{M}_2(\R^n)$ is an open problem, while Muriel Casalis has done two significant advances in it. The first advance concerns  the particular case $B$ and $C$ equal to zero. In that case  she has shown in \cite{CasaW} that the only solutions are the Wishart distributions on a symmetric cone.  Her second result  \cite{Casal} concerns  the particular case where $A(m)$ has rank one, namely has the form $amm^T.$ We call this subclass of  $\mathcal{M}_2(\R^n)$ the simple quadratic variance functions -or  equivalently the simple quadratic exponential families- and we denote it by $\mathcal{SM}_2(\R^n).$ We  speak also of the \textit{Casalis class.} A $G_0$-orbit is called a \textit{type} by Muriel Casalis.

\subsection{The $2n+4$ $ G_0$ orbits of the Casalis class}

 Let us describe now the  $G_0$-orbits of this Casalis class, by giving an element of each orbit attached to  its variance function, while the generating measures are given explicitely in Casalis \cite{Casal}.

\begin{itemize}
\item A first group is made of  $n$ $G_0$-orbits called Gaussian-Poisson, they are  such that $A(m)\equiv 0.$  For $k=1,2,\ldots,n$ their respective variances are

$$V_{I,k}(m)=\mathrm{diag}(m_1,m_2,\ldots, m_k, 1,\ldots,1)$$ with $M_F=(0,\infty)^k\times \R^{n-k}.$
\item The second  group  comprises $n$ orbits called Negative Multinomial - Gamma, denoted $(NM-ga)_k.$ They are such that,
for $k=1,2,\ldots,n,$ their respective variances are $$V_{IV,k}(m)=mm^T+\mathrm{diag}(0,m_2,\ldots, m_k, m_{1},\ldots,m_{1})$$ with $M_F=(0,\infty)^{k}\times \R^{n-k}.$

\item The third   group comprises four sporadic cases.
The first one is the Gaussian case $ V_{I,0}(m)=I_n.$

The second  one is the multinomial case, such that $A(m)=-mm^T.$
$$V_{II}(m)=-mm^T+\mathrm{diag}(m_1,m_2,\ldots, m_{n})$$ and $M_F$ is  the interior of the tetrahedron which is the convex hull of the set
$(0,e_1,\ldots,e_n)$ where $(e_1,\ldots,e_n)$ is the canonical basis of $\R^n.$

The third  one is the Negative multinomial case such that $A(m)=mm^T.$
$$V_{III}(m)=mm^T+\mathrm{diag}(m_1,m_2,\ldots, m_n)$$ with $M_F=(0,\infty)^{n}.$

The fourth  one is the hyperbolic case:
$$V_{V}(m)=mm^T+\mathrm{diag}(m_1,m_2,\ldots, m_{n-1} ,1+\sum_{k=1}^{n-1}m_k)$$ with $M_F=(0,\infty)^{n-1}\times \R.$
\end{itemize}
We are now in position to define the Morris Mora class $\mathcal{SM}_3(\R^n)$ as the image of  the Casalis class $\mathcal{SM}_2(\R^n)$ by the action of the
set $\{T_g; \ g \in G\}.$ The corresponding $V$ are matrices whose entries are polynomials in $m=(m_1,\ldots,m_n)$ of total degree $\leq 3$.
This point is not as obvious as it was for $n=1$ and we devote the next section to its proof.

\subsection{$T_g(V)$ is a  matrix polynomial if $V$ is in  the Casalis class} We need  a long preparation for the proof of Theorem 5.2.
If $g=\left[\begin{array}{cc}A&b\\c^T&d\end{array}\right]$ is in $G$ and if $V_F$ is a variance function, it is necessary that $M_g$ (see \eqref{UG}) is not empty to have $T_g(V_F)$ defined. For this reason we consider the set $$E_g=\{m; c^Tm+d>0\}$$ which will be not empty if either $c\neq 0$ or if $d>0$. In other terms, $E_g$  is empty if and only if $c=0$ and $d\leq 0.$   We denote by $\tilde {G}$ the set of $g\in G$ such that $E_g$ is not empty. In particular, $G_0\subset \tilde{G}.$  Here is the corresponding partition of $G:$

$$ \begin{array}{ccccc}&\rule{0.3mm}{6mm}&d>0&d=0&d<0\\----&
&----&----&----\\
c\neq 0&\rule{0.3mm}{6mm}&\tilde {G}&\tilde {G}&\tilde {G}\\
c= 0&\rule{0.3mm}{6mm}&G_0\subset \tilde {G}&\emptyset&G\setminus\tilde {G}
\end{array}$$

In the sequel we consider the following three elements of $G\setminus G_0$
\begin{eqnarray}\nonumber
g_c=\left[\begin{array}{cc}I_n&0\\c^T&1\end{array}\right],&& g_{c-}=\left[\begin{array}{cc}I_n&0\\c^T&-1\end{array}\right]\\ g_{b,c}=\left[\begin{array}{cc}I_n+bc^T&b\\c^T&1\end{array}\right]&=&
\left[\begin{array}{cc}I_n&b\\0&1\end{array}\right]\left[\begin{array}{cc}I_n&0\\c^T&1\end{array}\right]\label {GOODG}
\end{eqnarray}
and also the following three subsets of $G\setminus G_0$
\begin{eqnarray}\label {GOODH}
H_+&=&\{g_c; c\in \R^n,\  c\neq 0\},\\ H_-&=&\{g_{c-}; c\in \R^n,\ c\neq 0\},\nonumber\\H_0&=&\{g_{b,c}; b,c\in \R^n,\ c^Tb+1=0\}
\nonumber
\end{eqnarray}

Finally it is useful to introduce here the group
\begin{equation}\label{GROUPEH}
H=\{g=\left[\begin{array}{cc}I_n&0\\c^T&\lambda\end{array}\right]\, ; \lambda>0,\ c\in \R^n\}
\end{equation}
which will be essential in Section 7 below. Note that $$\left[\begin{array}{cc}I_n&0\\c^T&\lambda\end{array}\right]=     \left[\begin{array}{cc}I_n&0\\c^T&1\end{array}\right]\left[\begin{array}{cc}I_n&0\\0&\lambda\end{array}\right]=g_cJ_{\lambda}$$

In the sequel we will frequently use the simple formula for $h_{g_{b,c}}(m)=\frac{(I_n+bc^T)m+b}{c^Tm+1}:$

\begin{equation}\label{HGBC}  [h'_{g_{b,c}}(m)]^{-1}=(c^Tm+1)[I_n+mc^T]
\end{equation}

\vspace{4mm}\noindent\textbf{Proposition 5.1.}\begin{enumerate}
\item If $g\in G\setminus G_0$ then there exists $u$ and $v$ in $ \R^n$, and $g_0\in G_0$ such that $g=g_{v,u}g_0.$
\item $\tilde {G}=H_+ G_0\cup H_0 G_0\cup H_-G_0$
\item  If $c,b,b_0$ are in $\R^n$ such that $c^Tb=c^Tb_0$ there exists $g_0\in G_0$ such that $g_{c,b}=g_{c,b_0}g_0.$

\end{enumerate}

\vspace{4mm}\noindent\textbf{Proof.} Part 1) Let $g=\left[\begin{array}{cc}A&b\\c^T&d\end{array}\right]$ with $c\neq 0.$ We will show the existence of $u,v,A_1,b_1,d_1$ such that
$$\left[\begin{array}{cc}A&b\\c^T&d\end{array}\right]=\left[\begin{array}{cc}I_n+uv^T&u\\v^T&1\end{array}\right]\left[\begin{array}{cc}A_1&b_1\\0&d_1\end{array}\right],$$  in other terms we have to solve  the four equations
$$A_1+uv^TA_1=A,\ v^TA_1=c^T,\ b_1+uv^Tb_1+ud_1=b,\ v^Tb_1+d_1=d$$
We discuss three cases: $A=0$, $A$ invertible, and $A$ non zero but non invertible.

 We begin by a remark: if $g=\left[\begin{array}{cc}A&b\\c^T&d\end{array}\right]$ then the rank of $A$  is either $n-1$ or $n$. To see this  we consider
$$gg^T=\left[\begin{array}{cc}AA^T+bb^T&Ac+bd\\c^TA^T+db^T&\|c\|^2+d^3\end{array}\right]$$
Since $gg^T$ is positive definite, the rank of the principal matrix $AA^T+bb^T$ is $n$ which implies that the ranks of $AA^T$ and of $A$ cannot be $<n-1.$

\textsc{The case $A=0.$} It occurs only if $n=1$ from the above remark. Therefore , for $n=1$ the desired decomposition is

$$g=\left[\begin{array}{cc}0&b\\c^T&d\end{array}\right]=\left[\begin{array}{cc}0&b\\-1/b&1\end{array}\right]\times\left[\begin{array}{cc}-bc&b(1-d)\\0&1\end{array}\right].$$

\textsc{The case  $A$ is invertible}. Observe that $d-c^TA^{-1}b\neq 0$. This comes from a Cholesky decomposition of $g.$
Suppose now  $d-c^TA^{-1}b>0.$ Then $d_1=d-c^TA^{-1}b$, $b_1=b,$  $A_1=A  $, $u=0$  and $v=(A^{-1})^Tc$ will fit.

Suppose next that   $d-c^TA^{-1}b<0.$ Suppose first that $d\neq 0$. We chose $d_1=d$, $b_1=0$ and  $u=b/d. $  This gives $A_1=A-\frac{1}{d}bc^T$ which is invertible, again by another  Cholesky decomposition. Therefore finally we take $v^T=c^TA_1^{-1}.$

Consider the case $d=0$. Since $ b_1+u(v^Tb_1+d_1)=b,\ v^Tb_1+d_1=0$ one has $b_1=b$ and we have to chose $A_1,u,v,d_1>0$ such that
$$A_1+uv^TA_1=A,\ v^TA_1=c^T, \ v^Tb+d_1=0.$$ A delicate point is the fact that $c^TA^{-1}b>0$ and we want $d_1>0.$ We chose $u=-\lambda Ac$ such that $\lambda c^Tc>1.$ Therefore $A_1=A-uc^T=A(I-\lambda cc^T)$ and
$$-d_1=v^Tb=c^T(I_n-\lambda cc^T)^{-1}A^{-1}b =c^T(I+\frac{\lambda cc^T}{1-\lambda c^Tc})A^{-1}b=c^TA^{-1}b\times \frac{1}{1-\lambda c^Tc}<0.$$

 \textsc{The case $A$ is non zero and non invertible}. We introduce $$g^{-1}=\left[\begin{array}{cc}A'&b'\\(c')^T&d'\end{array}\right].$$
The fact that $AA'+b(c')^T=I_n$ implies that none of $A',b,c'$ is zero. Let $v\in \R^n\setminus \{0\} $ such that $v^TA=0.$ Since $v^Tb(c')^T=v^T$
then $v^Tb\neq 0$ and we can chose $v$  such that  $v^Tb=-1.$  If  $d\neq 0,$
take $u=b/d$ , $A_1=A-\frac{1}{d}bc^T,$ $b_1=0$ and $d_1=d$.
Again, $A_1$ is invertible by a Cholesky decomposition.

Finally we are left with the case where $d=0$ (and $A$ non invertible). This is the last case that we have to solve for ending the proof of part 1. As seen above $A$ has rank $n-1$. To this end we assume first that $A$ is diagonal of the form
$$D=\mathrm {diag}(d_1,\ldots,d_{n-1},0)=\mathrm{diag}(D_1.0)$$ where $D_1$ is invertible. We want  to find a $A_1$ invertible, $  u,v,b_1,d_1>0$ such that

$$g=\left[\begin{array}{cc}D&b\\c^T&0\end{array}\right]=\left[\begin{array}{cc}I_n+uv^T&u\\v^T&1\end{array}\right]\times\left[\begin{array}{cc}A_1&b_1\\0&d_1\end{array}\right]$$ or
\begin{equation}\label{OUF}D=A_1+uv^TA_1,\ c^T-v^TA_1,\ b=b_1+u(v^Tb_1+d_1),  \ v^Tb_1+d_1=0.\end{equation}
We get necessarily  $b=b_1.$ Let us also observe that  if $c_n=0$ then the column $n$ of the matrix $g$ is zero, which contradicts invertibility of $g.$ Similarly $b_n=0$ produces a row of zeros in $g$, another impossibility. Thus $b_nc_n\neq 0.$ Now for filling the conditions of \eqref{OUF} we choose $u^T=(0,\ldots,0,\lambda).$ Also for simplicity we write $c^T=(C^T,c_n)$ and $b^T=(B,b_n)$
From \eqref{OUF} we get $A_1=D-uc^T.$ Let us show that with the above choice of $u$ then $A_1$ is invertible and let us compute $A_1^{-1}.$
With the above notations we have

$$A_1=\left[\begin{array}{cc}D_1&0\\\lambda C^T&\lambda c_n\end{array}\right], \ \ A_1^{-1}=\left[\begin{array}{cc}D_1^{-1}&0\\-\frac{1}{c_n}C^TD_1^{-1}&\frac{1}{\lambda c_n}\end{array}\right].$$
We are now able to compute $v$ corresponding to this choice of $u$ and we get
$$0=v^Tb+b+d_1=c^TA_1^{-1}b+d_1=K+d_1+\frac{b_n}{\lambda}$$  where $K$ is a number the  value of which does not import. Finally from the last equality we can choose  a $\lambda$  such that $d_1>0.$

We now come to the general case where $A$ is not necessarily diagonal but has rank $n-1.$ We use  a polar decomposition of $A=UDV^T$ where $U$ and $V$ are orthogonal matrices and $D$ has the form above.

\begin{eqnarray*}\left[\begin{array}{cc}A&b\\c^T&0\end{array}\right]&=&\left[\begin{array}{cc}UDV^T&UU^Tb\\c^TV^TV&0\end{array}\right]=
\left[\begin{array}{cc}U&0\\0&1\end{array}\right]\left[\begin{array}{cc}D&U^Tb\\(V^Tc)^T&0\end{array}\right]\left[\begin{array}{cc}V^T&0\\0&1\end{array}\right]\\&=&
\left[\begin{array}{cc}U&0\\0&1\end{array}\right]\left[\begin{array}{cc}U^T&0\\0&1\end{array}\right]\left[\begin{array}{cc}U&0\\0&1\end{array}\right]\left[\begin{array}{cc}A_1&b_1\\0&d_1\end{array}\right]\left[\begin{array}{cc}V^T&0\\0&1\end{array}\right]\\&=&
\left[\begin{array}{cc}I_n+(Uu)(Vv)^T&Uu\\(Vv)^T&1\end{array}\right]\times \left[\begin{array}{cc}A_1V^T&b_1\\0&d_1\end{array}\right]
\end{eqnarray*}

Part 2) If  $g\in G_0$  then $g=I_nG\in H_{1}G_0.$ Consider $g=\left[\begin{array}{cc}A&b\\c^T&d\end{array}\right]$ with $c\neq 0.$

If $g\in \tilde{G}\setminus G_0$ , that is if $c\neq 0$ then from part 1 there exist $u,v\in R^n$ and $g_0$ such that

$$g=\left[\begin{array}{cc}I_n+uv^T&u\\v^T&1\end{array}\right]g_0=g_{u,v}g_0.$$  If $v^Tu+1=0$  then $g_{u,v}\in H_0.$
If  $v^Tu+1=R\epsilon$ with $\epsilon =\pm 1$ and $R>0$ we write
$$ g_{u,v}=\left[\begin{array}{cc}I_n+uv^T&u\\v^T&1\end{array}\right]=\left[\begin{array}{cc}I_n&0\\v^T/R&\epsilon\end{array}\right]\left[\begin{array}{cc}I_n+uv^T&u\\0&1/R\end{array}\right].$$Since the eigenvalues of $I_n+uv^T$ are $v^Tu+1\neq 0$ and $1$ with multiplicity $n-1$ then  $I_n+uv^T$  is invertible  and  $\left[\begin{array}{cc}I_n+uv^T&u\\0&1/R\end{array}\right] $ belongs to $G_0.$
Since $\left[\begin{array}{cc}I_n&0\\v^T/R&\epsilon\end{array}\right]$ is in $H_{\epsilon}$ we have shown $\tilde {G}\subset H_+ G_0\cup H_0 G_0\cup H_+ G_0.$ The other inclusion  is clear.

Part 3) This comes from the equality $$ \left[\begin{array}{cc}I_n+bc^T&b\\c^T&1\end{array}\right]=\left[\begin{array}{cc}I_n+b_0c^T&b_0\\c^T&1\end{array}\right]\left[\begin{array}{cc}I_n+(b-b_0)c^T&b-b_0\\0&1\end{array}\right]$$

\vspace{4mm}\noindent\textbf{Theorem 5.2.} If $g\in G$ and $F\in \mathcal{SM}_2(R^n)$ then $T_g(V_F)$ is a matrix polynomial of degree $\leq 3. $

\vspace{4mm}\noindent\textbf{Proof.} From part 1) of Proposition 5.1, any $g$ in $G$ can be written as

$$g=\left[\begin{array}{cc}I_n&b\\0&1\end{array}\right]\left[\begin{array}{cc}I_n&0\\c^T&1\end{array}\right]g_0.$$
 where $g_0$ is in $G_0$. Since the elements of $G_0$ preserve the degrees, and since from Proposition 3.2 we have $T_g\circ T_{g_1}=T_{g_1g}$
enough is to prove the proposition for $g=\left[\begin{array}{cc}I_n&0\\c^T&1\end{array}\right],$ with the help of  Proposition 3.1:
\begin{eqnarray*}h_g(m)=\frac{m}{c^Tm+1}&,&\ \ h'_g(m)=\frac{1}{(c^Tm+1)^2}[(c^Tm+1)I_n-mc^T],
\\ h'_g(m)^{-1}=(c^Tm+1)[I_n+mc^T]&,&\  \ ( h'_g(m)^T)^{-1}=(c^Tm+1)[I_n+cm^T]
\end{eqnarray*}
Furthermore, from definition \eqref{TG} the map $V\mapsto T_g(V)$ is linear, enough is to  prove the result for $A(m)=m^Tm, $ for $B(m)=B_1m_1+\cdots+B_nm_n$ and for the constant matrix $C.$

 Calculation gives  the remarkable formulas  \begin{eqnarray}\label{Q2}V(m)=mm^T &\Rightarrow &T_g(V)(m)=(c^Tm+1)mm^T\\ \label{Q1} V(m)=B(m)&\Rightarrow&T_g(V)(m)=[I_n+mc^T]B(m)[I_n+cm^T]\\\label{Q0}
V(m)=C &\Rightarrow & T_g(V)(m)=(c^Tm+1)[I_n+mc^T]C[I_n+cm^T].\end{eqnarray}  In the three cases we obtain polynomials in  $m$ of degree $\leq 3.$

\vspace{4mm}\noindent\textbf{Comments.} 1) It worths to mention  that if we apply the above $T_g$ to the variance function $V(m)=\lambda P(m)$ of the family of Wishart distributions (here $P(m)$ is the linear operator defined on the linear space $S$ of symmetric matrices  of order $n$ by $\Delta\mapsto P(m)(\Delta)=m\Delta m$ ) then we have
\begin{eqnarray*}T_{g}(V)(\Delta)&=& \frac{\lambda}{c^Tm+1}(\mathrm{id}_S+m\otimes c)m(\mathrm{id}_S+c\otimes m)\Delta m\\
&=&\frac{\lambda}{c^Tm+1}\times (P(m)+mc\otimes m^2+m\otimes (P(m)c)+( m^Tc^2) \, m\otimes m^2)(\Delta)
\end{eqnarray*}
where  for $c,m\in S$, the symbol $c^Tm$ is the trace of $cm$ and the  symbol $m\otimes c$ is the linear map of $S$ into $S$ defined by $\Delta\mapsto mc^T\Delta.$
The interesting point is that $T_{g}(V)$ is not a matrix polynomial anymore.

2) A similar remark holds if we want to apply $T_g$ to a diagonal variance of the form $V(m)=\mathrm{diag}(V_1(m_1),\ldots,V_n(m_n))$  where the $V_i(m_i)=am_i^2+bm_1+c$ belong to the Morris class with $a\neq 0.$ The result will be a matrix with rational, not polynomial coefficients.

\subsection {The $n+3$ $G$-orbits of the Morris Mora class in $\R^n$}
Recall that  the Morris Mora class in $\R^n$ is the set of images by $G$ of the Casalis class and, from Theorem 5.2, all their variances are matrix polynomials of degree $\leq 3.$ Since $G_0\subset G$  the technique for describing the $G$-orbits of the Morris Mora is to consider each of the $2n+4$ $G_0$-orbits of the Casalis class and to compute its image by $G$. As we have seen in the case $n=1$, the number of $G$-orbits in the Morris Mora class is less than the number of  $G_0$-orbits in the Casalis class. This subsection  will show in Theorem 5.4  that for fixed $k\in \{1,\ldots,n\} $, $V_{I,k}$ and $V_{IV,k}$  are in the same $G $- orbit, that $V_{II}$ and $V_{III}$  are in the same orbit  and  finally that  the $n+3$ $G$-orbits  of $V_{I,0},\ldots,V_{I,n}, V_{II}, V_{V}$ are disjoint.  May be it is good to recall here that for $n=1$ we have already seen that in the Morris class, informally $m$ and $m^2$ were is the same $G$-orbit, that $-m^2+m $ and $m^2+m$ were in the same $G$-orbit while $1$ and $m^2+1$ were alone.

 The rest of the task  is devoted to the description of these cubic new  exponential families of the set $G(\mathcal{SM}_2)\setminus \mathcal{SM}_2$ for each of these $n+3$ $G$-orbits. This work  will be performed in Part two. Several tools will have still to be gathered in Sections 6 and 7.

We begin with a not very appealing technical proposition.

\vspace{4mm}\noindent\textbf{Proposition 5.3.} Let $g_{b,c}=\left[\begin{array}{cc}I_n+bc^T&b\\c^T&1\end{array}\right]$ as in \eqref{GOODG} and $A$ the square matrix of order $n$ written by blocks as $A=\mathrm{diag}(I_{i-1}, J,I_{n-j})$ where $ (x_i,\ldots,x_j)J=(x_j,x_{i+1},\ldots,x_{j-1},x_i)$  and $i<j.$ In other terms, $A$ is the permutation matrix of the transposition $(i,j)$. Let $g_0=\mathrm{diag}(A,1).$
Denote also $ D(m)=\mathrm{diag}(m_1,\ldots,m_n) $ and $E(m)=mm^T.$ Then for all $i<j$
\begin{eqnarray}
T_{g_0}(T_{g_{b,c}})(D)(m)&=& T_{g_{Ab,Ac}}(D)(m)\\
T_{g_0}(T_{g_{b,c}})(E)(m)&=& T_{g_{Ab,Ac}}(E)(m)
\end{eqnarray}

\vspace{4mm}\noindent\textbf{Proof.} Proposition 3.2 implies that $T_{g_0}(T_{g_{b,c}}(V))=T_{g_1}(V)$ with  $ g_1=\left[\begin{array}{cc}A+bc^TA&b\\c^TA&1\end{array}\right].$
while $g_{Ab,Ac}=\left[\begin{array}{cc}I_n+Abc^TA&Ab\\c^TA&1\end{array}\right].$\\
We get by careful computation  using \eqref{HGBC} \begin{eqnarray}h_{g_1}(m)=\frac{(I_n+bc^T)Am+b}{c^TAm+1} \ \  \ \  , \ \  ( h'_{g_1})^{-1}(m)=(c^TAm+1)(A+mc^T)\ \  \ \   \ \  \\
h_{g_{Ab,Ac}}(m)=\frac{(I_n+Abc^TA)m+Ab}{c^TAm+1}\ \  , \ \ ( h'_{g_{Ab,Ac}})^{-1}(m) =(c^TAm+1)(A+mc^T)A\end{eqnarray}
For simplification we introduce the scalar $k=c^TAm+1.$ We get
\begin{eqnarray*}T_{g_1}(D)(m)&=&\frac{1}{k}\times k\, (A+mc^T)D\left(\frac{1}{k}((I_n+bc^T)Am+b)\right)(A+cm^T)k\\&=&(A+mc^T)D((I_n+bc^T)Am+b))(A+cm^T)\end{eqnarray*}
$$T_{g_{Ab,Ac}}(D)(m)=(A+mc^T)AD((I_n+Abc^TA)Am+Ab))A(A+cm^T)$$
Using the fact that $AD(Am)A=D(m)$ we get the equality of the last two lines.

\vspace{4mm}\noindent
  Similarly for $E(mm^T)$ we get in the same way that $T_{g_1}(E)(m)=T_{g_{Ab,Ac}}(E)(m).$

Here comes  the theorem of the $n+3$ classes for $G.$

\vspace{4mm}\noindent\textbf{Theorem 5.4.} For all $k=1,\ldots,n$  the two classes  $V_{I,k}$  and $V_{IV,k}$ are in the same $G$-orbit. Similarly the two classes $V_{II}$ and $V_{III}$ are in the same $ G$ orbit.  Finally the $ n+3$ $G$-orbits of $V_{I,0},\ V_{I,1}, \dots, \ V_{I,n},\ V_{II}$ and $V_{V}$ are disjoint.

\vspace{4mm}\noindent\textbf{Proof.}
From Proposition 5.1, the technique for describing the elements of  a $G$-orbit of  $V$ which are of the form  $amm^T+B(m)+C$ is to consider its image only  by all  $T_{g_{(b,c)}}$ namely
$$T_{g_{(b,c)}}(V)=(c^Tm+1)amm^T+(I_n+mc^T)B(m)(I_n+cm^T)+(c^Tm+1)(I_n+mc^T)C(I_n+cm^T).$$
But since in the present  theorem we are interested in the elements of the $G$-orbit of $V$ with no element of degree 3, we observe that $T_{g_{(b,c)}}(V)$ has  the form $a_1mm^T+B_1(m)+C_1$ if and only if  $V(m)=amm^T+B(m)+C$ is such that  for all $m$
\begin{equation}\label{NO3}
ac^Tm(mm^T)+mc^TB(m)\, cm^T+c^Tm(mc^T)C\, cm^T=0\end{equation}
which is trivially fulfilled if $c=0.$  In the sequel of the proof we shall use only $T_{g_{(b,c)}}$ satisfying \eqref{NO3}.

\vspace{4mm}\noindent\textsc{The Gaussian case $V_{I,0}$}. This is the Gaussian case since a  representative  is $V(m)=I_n.$  Clearly condition \eqref{NO3} is fulfilled if and only if $c=0$ . This shows that the only elements of the $G$-orbit of $V_{I,0}$ of the form $amm^T+B(m)+C$ are the constants, thus belonging to the $G_0$-orbit of $V_{I,0}.$

\vspace{4mm}\noindent\textsc{The case Multinomial and Negative Multinomial $(V_{II},\ V_{III})$}. For considering the action of $\tilde{G}$ to the cases $V_{II}(m)=-E(m)+D(m)$ and $V_{II}(m)=E(m)+D(m)$ we use the decomposition of Proposition 5.1, part 2, by considering  only the transforms $T_{g_c+}$, $T_{g_c-}$ and $T_{g_{b,c}}$ applied to $-E(m)+D(m)$  and wondering which ones have the form
$a_1mm^T+B_1(m)+C_1.$ Using the formulas \eqref{Q2} and \eqref{Q1} we have
\begin{equation}\label{ED}T_g(E)(m)=(c^Tm+1) E(m), \ T_g(D)(m)=(I_n+mc^T) D(m)(I_n+cm^T)\end{equation}
Therefore $T_g(-E+D)$ is quadratic if and only if $-c^TmE(m)+mc^TD(m)cm^T=0$ or more explicitely
$$(-(c_1m_1+\cdots +c_nm_n)+(c_1^2m_1+\cdots+c_n^2m_n))E(m)=0$$ that is to say $c_i=0$ or $1$ for $i=1,\ldots,n.$
Note that from $\eqref{ED}$ if $c_i=1$ for all $i=1,\ldots,n$ we have $T_g(-E+D)=E+D$ or $T_g(V_{II})=T_g(V_{III}).$
Since here $c\neq 0$ the set of $i$ such that $c_i=1$ is not empty. Without loss of generality we can assume that $c=(1,1, \ldots,1,0,\ldots,0)$ containing $k>0$ ones and $n-k$ zeros.

We now show that $k=n.$ If not let us concentrate on the (2,2) submatrix of $T_g(-E+D)$ corresponding to lines and columns $k$ and $k+1$. From $\eqref{ED}$ and a simple calculation this submatrix is
$$\left[\begin{array}{cc}m^2_k+m_k&m_k\, m_{k+1}\\m_k\, m_{k+1}&-m^2_{k+1}+m_{k+1}\end{array}\right].$$
The fact that the coefficients of $m^2_{k+1}$ and $m^2_k$ are opposite prevents $T_g(-E+D)$  to be a simple quadratic variance. Thus $k=n$. In particular $V_{II}$ and $V_{III}$ are in the same $G$-orbit, with none other members of the set of the $G_0$-orbits of the Casalis class.

\vspace{4mm}\noindent\textsc{The  Poisson -Gaussian and NG-gamma $(V_{I,k},\ V_{IV,k})$ cases.}  In this part, where $k$ is fixed in $\{1,\ldots,n\}$ we write for $x\in \R^n$

$$x=\left(\begin{array}{c}x_-\\ x_+\end{array}\right)\in \left(\begin{array}{c}\R^k\\ \R^{n-k}\end{array}\right).$$

Taking $V(m)=\mathrm{diag }(m_1,\ldots,m_k,1,\ldots,1)$ condition \eqref{NO3} for applying $T_{g_{(b,c)}}$ to $V$ reads

\begin{eqnarray}\label{HA1}
(c_-^T\mathrm{diag}( b_- )c_-+c_+^Tc_+)c_i+ c_i^2&=&0\ \ \mathrm{for}\ \ 1\leq i \leq k;\\
\label{HA2}
(c_-^T\mathrm{diag}( b_- )c_-+c_+^Tc_+)c_i&=&0\ \ \mathrm{for}\ \ k+1\leq i \leq n.
\end{eqnarray}
Let show that \eqref{HA1} and \eqref{HA2} imply $c_+=0.$ If $c_+\neq 0$ then there exists $i_0>k$ such that $c_{i_0}\neq 0$ then
$c_-^T\mathrm{diag}( b_- )c_-+c_+^Tc_+=0$, then $c_i^2=0$ for all $i\leq k,$  then $c_+^Tc_+=0$ a contradiction. Since we assume $c\neq 0$
let $i_1\leq k$ such that $c_{i_1}\neq 0.$ From Proposition 4.1 part 3) we can choose $b$  such that $b_i\neq 0$ if and only if $i=i_1.$ From \eqref{HA1} we have $(c_{i_1}^2b_{i_1}+c_i)c_i=0$ for all $i\leq k.$ In particular $1+c_{i_1}b_{i_1}=0$. Therefore $g_{b,c}$ is in $H_0$ (Recall that $H_0$ is defined in \eqref{GOODH}). From Proposition 5.3 without loss of generality we may assume that $i_1=1,$ that $b_1=-1/c_1$ and that
the $k_0\leq k$ non zero elements of $c_1,\ldots,c_k$ are the first. Remark also that $c_1<0$ to insure that $(M_F)_{g_{b,c}}$ is not empty.
 Under these circumstances  we write by blocks
$$T_{g_{b,c}}(V_{I,k})=\left[
\begin{array}{cc}A_0(m)&B_0(m)\\B_0(m)^T&C_0(m)\end{array}\right]$$ where the symmetric matrix $A_0(m)$ of order $k$  and $B_0$ and $C_0$ are
\begin{eqnarray}
A_0(m)&=&\nonumber|c_1|mm^T+\mathrm{diag}(m_-)\, c_-m_-^T+m_-c_-^T\, \mathrm{diag}(m_-)\\&&\label{AZERO}+(1+c_-^Tm_-)(e_1m_-^T+m_-e_1^T-b_1\mathrm{diag (e_1)})+\mathrm{diag}(m)\\\nonumber
&=&|c_1|\left[\begin{array}{ccc}m_1^2&\ldots&m_1m_k \\\ldots&\ldots&\ldots\\m_km_1&\ldots&m_k^2\end{array}\right]+\left[\begin{array}{ccc}2c_1m_1^2&\ldots&(c_1+c_k) m_1m_k\\\ldots&\ldots&\ldots\\(c_k+c_1) m_km_1&\ldots&2c_km_k^2\end{array}\right]\\&&-(1+c_-^Tm_-)\left[\begin{array}{cccc}2m_1-b_1&m_2&\ldots&m_k \\m_2&0&\ldots&0\\\ldots&\ldots&\ldots&\ldots \\m_k&0&\ldots&0\end{array}\right]+\left[\begin{array}{ccc}m_1&\ldots&0 \\\ldots&\ldots&\ldots\\0&\ldots&m_k\end{array}\right],\nonumber
\end{eqnarray}
\begin{equation}\label{BZERO}B_0(m)=|c_1|m_-m_+^T+\mathrm{diag}(m_-)m_-c_+^T-(1+c_-^Tm_-)e_1m_+^T,\end{equation}
\begin{equation}\label{CZERO}C_0(m)=|c_1|m_+m_+^T+(1+c_-^Tm_-)I_{n-k},\end{equation}
From \eqref{AZERO}, \eqref{BZERO} and \eqref{CZERO} one sees that $T_{g_{b,c}}(V_{I,k})$ is simple quadratic if and only if $c_2=c_3=\ldots=c_k=0.$ For simplification,  we consider the vector $M^T=(m_2,\ldots,m_k)$ and the number  $K=1+c_-^Tm_-$. With this notation we have now, splitting the variance matrix with the partition $(1, k-1,n-k)$
\begin{equation}\label{BIGONE}
T_{g_{b,c}}(V_{I,k})(m)=\left[\begin{array}{ccc}
|c_1|m_1^2-2m_1+b_1&-KM^T&-Km_+^T\\-KM&|c_1|MM^T+\mathrm{diag(M)}&|c_1|Mm_+^T\\-Km_+&|c_1|m_+M^T&|c_1|m_+m_+^T+KI_{n-k}
\end{array}\right]\end{equation}
Applying the transforms $$g_1=\left[\begin{array}{cc}I_n&0\\0&|c_1|\end{array}\right],\  g_2=\left[\begin{array}{cc}I_n&e_1\\0&1\end{array}\right], \ g_3=\left[\begin{array}{cc}\frac{1}{|c_1|}\mathrm{diag}(e_1)&0\\0&1\end{array}\right]$$ we get

\begin{equation}\label{BIGTWO}
T_{g_1}(T_{g_{b,c}}(V_{I,k}))(m)=\left[\begin{array}{ccc}
(m_1-1)^2&(m_1-1)M^T&(m_1-1)m_+^T\\(m_1-1)M&MM^T+\mathrm{diag(M)}&Mm_+^T\\(m_1-1)m_+&m_+M^T&m_+m_+^T+c_1(m_1-1)I_{n-k}
\end{array}\right]\end{equation}

\begin{equation}\label{BIGTHREE}
T_{g_2}(T_{g_1}(T_{g_{b,c}}(V_{I,k})))(m)=mm^T+\mathrm{diag}(0,M,c_1m_1I_{n-k})
\end{equation}
which is the desired result.

\vspace{4mm} \textsc{The hyperbolic case $V_{V}$}. Here there is no calculation to do since  $V_{V}$ never appears in the other $n+2$ orbits.

\section{The natural exponential families concentrated on $\N^n$}
Thus, our ultimate aim is  to the description of the elements of the $n+3$ $G$-orbits of the Casalis class in Part two.  Many of them have representatives which are concentrated on $\N^n$ (where $\N=\{0,1,2,\ldots\}$). Note that the fact that an exponential family is concentrated on $\N^n$ does not necessary implies that all the elements of its $G$-orbit are concentrated on a discrete set: recall that the $G$-orbit of the one dimensional Poisson family contains the Gamma families. However, the classification that we are performing  makes a large use of the NEF concentrated on $\N^n.$

This section uses Lagrange formula in several dimensions for giving a characterization among the variance functions $V_F$  of the ones corresponding to $F$ concentrated on $\N^n$ (Proposition 6.3).

\subsection{The Lagrange formula. }We first adapt the description of the classical objects attached to a natural exponential family $F$ when it is concentrated on $\N^n.$  In this case, the information on a generating measure $(\mu(k))_{k\in \N^n}$ of $F$ is generally given by its generating function
$$f_{\mu}(z)=\sum _{k\in \N^n}\mu(k)z^k$$ with the notation $z^k=z_1^{k_1}\ldots z_n^{k_n}.$ In many circumstances in this section  we will have to consider functions of the form
$$f_{\mu}(w)=\sum _{k\in \N^n}w^k[z^k](g_k(z)).$$
where $[z^k](g(z))$ means the coefficient of $z^k=z_1^{k_1}\ldots z_n^{k_n}$ in the power series expansion of a analytic function $g$ on $\C^n$ or $\R^n$ defined on  a neiborhood of $0$.

We will need the notations
\begin{eqnarray}\label{JACOBIAN1}M(g)(z)&=&M(g_1,\ldots,g_n)(z)=I_n-\left[\frac{z_i}{g_i(z)}\times \frac{\partial g_i}{\partial z_j}(z)\right]_{1\leq i,j\leq n},\\\ D(g)(z)&=&\det(M(g)(z))\label{JACOBIAN2}\end{eqnarray}
This matrix will appear in the calculation of the derivative of
\begin{equation}\label{JACOBIAN3}z=(z_1,\ldots,z_n)\mapsto \left(\frac{z_1}{g_1(z)},\ldots,\frac{z_n}{g_n(z)}\right)\end{equation}
which is the matrix $$\mathrm{diag}  \left(\frac{z_1}{g_1(z)},\ldots,\frac{z_n}{g_n(z)}\right)        M(g_1,\ldots,g_n)(z)$$
When the transformation \eqref{JACOBIAN3} will be used for a change of variable the Jacobian will be
$$\left(\prod_{i=1}^n\frac{z_i}{g_i(z)}\right)  D(g_1,\ldots,g_n)(z).$$

We recall without proof the Lagrange theorem in Proposition 6.1. Our references are Good \cite{Good1} and Gessel \cite{Gesse}. Also the formal proof of  Goulden and Jackson \cite{Gould} page 21 is relatively easy to follow.

\vspace{4mm}\noindent\textbf{Proposition 6.1.} Let  $U$ be a connected open set of $\R^n$ containing 0 and $g=(g_1,\ldots,g_n)$ analytic on $U$ valued in $\C^n.$ We assume that for all $i=1,\ldots,n$ we have $ g_i(0)\neq 0.$ Under these circumstances there exists a neighborhood $W$ of $0$ in $\C^n$  and an analytic function $h$ from $W$ to $U$ such that
\begin{equation}\label{LAGRANGE1}h(w)=\mathrm{diag}(w)g(h(w))\end{equation}
meaning for all $i=1,\ldots$ $$h_i(w_1,\ldots,w_n)=h_i(w)=w_ig_i(h_1(w),\ldots, h_n(w)).$$
Furthermore, if $g_0: U\rightarrow \C$ is analytic  and  if for $z\in U$ the determinant $D(g_1,\ldots,g_n)(z)$  is not zero, then using the notation \eqref{JACOBIAN2} we have the Lagrange formula:
\begin{equation}\label{LAGRANGE3}
[w^k]g_0(h(w))=[z^k](g_0(z)g(z)^kD(g_1,\ldots,g_n)(z)).\end{equation}
In other terms
\begin{equation}\label{LAGRANGE4}
g_0(h(w))=\sum_{k\in \N^n}w^k[z^k](g_0(z)g^k(z)D(g_1,\ldots,g_n)(z).\end{equation}

\subsection{Variance functions on $\N^n$}

 In this subsection, Theorem 6.2 gives the important characterization of the variance fonctions of exponential families concentrated on $\N^n$ (with  the usual restriction $\mu_{e_i}\neq 0$ for all $i=1,\ldots,n$ and $\mu_{0}\neq 0.$ This theorem is the extension of the one dimensional result in Letac and Mora (1990) saying that  an exponential family $F(\mu)$ is concentrated on $\N$ such that $\mu(0)\neq 0$ and $\mu(1)\neq 0$ if and only if its variance function $V(m)$ is the sum of a power series with positive radius of convergence such that $V(m)=m+\sum_{k\geq 2}a_km^k.$ Through the Lagrange formula, Theorem 6.2 also gives the values of the $\mu_k.$

\vspace{4mm}\noindent\textbf{Theorem 6.2.} Let $F=F(\mu)$ be a NEF on $\R^n$ with variance function $V_F$ defined on $M_F.$ Then $F$ is concentrated on $\N^n$ with $\mu=\sum_{k\in \N^n}\mu_k\delta_k$ such that
$\mu_0,
\mu_{e_1},\ldots,\mu_{e_n}$ are positive if and only if the following two conditions are fullfilled:

\begin{enumerate}\item $M_F\subset (0,\infty)^n$ and $0\in\overline{M_F},$
\item There exists a simply connected open set $\mathcal{O}$ containing $M_F$ and $0$, and a real-analytic function
$ \phi'$ from $\mathcal{O}$ to $\R^n$ such that $ \phi'(0)=(1,\ldots,1)^T$ and for $m\in M_F,$ we have $$\phi'(m)=(V_F(m))^{-1}(m).$$\end{enumerate}

Under these circumstances $\phi'$ is the differential of a real-analytic function $\phi$ from $\mathcal{O}$ to $\R.$ Furthermore there exists a neighboorhood $\mathcal{V}$ of $0$ in $\C^n$ and an analytic function $G=(G_1,\ldots,G_n)$ from $\mathcal{V}$ to $\C^n$ such that $G(0)=(1,\ldots,1)^T$ and \begin{equation}\label{FUNCTIONG}
\mathrm{diag}\left(\frac{1}{G_1(m)},\ldots,\frac{1}{G_n(m)}\right)\times G'(m)(m)=(1,\ldots,1)^T-\phi'(m).
\end{equation}
Finally $\mu_0=e^{\phi(0)}$ and using the notation \eqref{JACOBIAN2}:

\begin{equation}\label{LAGRANGEMU}\mu_{k}=[z^k]e^{\phi(z)}G(z)^kD(G(z))\end{equation}

 Theorem 6.2 will be the consequence of Propositions 6.3, 6.4, 6.5 below.

\vspace{4mm}\noindent\textbf{Proposition 6.3.} Let $\mu(dx)=\sum_{k\in \N^n}\mu_k\delta_{k}(dx)$  a positive measure of $\N^n$ such that $\mu_0,\mu_{e_1},\ldots,\mu_{e_n}$  are all $>0$ and such that the domain of convergence $D\subset \C^n$ of $f(|z|)=\sum _{k\in \N^n}\mu_{k}|z|^k$ contains $0$ in its interior $D^{o}.$ Call $F$ the natural exponential family generated by $\mu.$ Then
\begin{enumerate}
\item $M_F\subset (0,\infty)^n$ and $0\in \overline{M_F}.$
Furthermore there exists an open set $U \subset D^o$  containing $D^o\cap (0,\infty)^n$ such that for all $z\in U$ and for all $i=1,\ldots,n$ one has $\frac{\partial f}{\partial z_i}(z)\neq 0.$

\item If $g: U \rightarrow \C^n$ is defined by  $g_i(z)=f(z)/\frac{\partial f}{\partial z_i}(z)$ then there exists an open subset $\mathcal{O}\subset \C^n$ simply connected containing $0$ and $M_F$, as well as a function $h: \mathcal{O}\rightarrow  U$ such that for all $m$ in $\mathcal{O}$ one has
\begin{equation}\label{LAGRANGE}h(m)=\mathrm{diag}(m)g(h(m)).\end{equation}

\item With these notations one has for all $m\in M_F$
\begin{eqnarray}\nonumber
P(m,F)(dx)&=&\sum _{k\in \N^n}\frac{h(m)^k}{f(h(m))}\mu_{k}\delta_k(dx)\\
V_F(m)&=&(h'(m))^{-1}\mathrm{diag}(h(m)).\label{VANDh}
\end{eqnarray}
\end{enumerate}
\vspace{4mm}\noindent\textbf{Proof.}

\vspace{4mm}\textsc{Proof of part 1.}
Consider $$U=\left\{z\in D^o; \frac{\partial f}{\partial z_i}(z)\neq 0 \ \ \forall\   i=1,\ldots,n \right\}.$$ Since
$$f(z)=\mu_0+\mu_{e_1}z_1+\cdots+\mu_{e_n} z_n+\sum_{k\in \N^n,\ |k|\geq 2}\mu_kz^k$$ the fact that $\mu_{e_i}>0$ implies $0\in U.$
The set $U$  is open and such that $D^o\cap(0,\infty)^n\subset U\subset D^o.$

In the sequel, for simplification, we rather write $L=L_{\mu}$, $k=k_{\mu}$, $\psi=\psi_{\mu}.$
Denoting $e^{\theta}=(e^{\theta_1},\ldots,e^{\theta_n})$  the Laplace transform of $\mu$ is
$$L(\theta)=f(e^{\theta})=\sum _{k\in \N^n}\mu_{k}e^{\<k,\theta\>}.$$ Hence the cumulant function is $k(\theta)=\log L(\theta)$ and $$\frac{\partial k}{\partial \theta_i}=m_i=\frac{e^{\theta_i}}{g_i(e^{\theta})}.$$ Since there exists $R>0$ such that
$\{\theta; \ \|e^{\theta}\|<R\}$ is contained in $ \Theta(\mu)$  then doing $\theta_1=\ldots=\theta_n\to -\infty$ shows that $0\in \overline {M_F}$ since $\mu_0\neq 0$.

\vspace{4mm}\textsc{A fact of topology for part 2.}   It says  that if $U_1$ and $B$ are two simply connected open sets of $ \R^n$ with a connected non empty intersection, then $U_1\cup B$ is also simply connected. To see this, let $\Gamma$
be a loop contained in $U_1\cup B$ such that $\Gamma\cap U_1$ and $\Gamma\cap B$ are not empty, and let $C$ be a connected component of $\Gamma\cap B.$ Let $A_1$ and $A_2$ the extreme points of $C$ : they belong to $\overline {B}\setminus B.$
We have also $A_1$ and $A_2$ in $U_1.$

Let $B_1$ and $B_2$ two balls centered in $A_1$ and $A_2$ respectively both contained in $U_1$. Let $P_1\in B_1\cap C$ and
$P_2\in B_2\cap C\subset U_1\cap B$. Let $\gamma_C$ a path linking $P_1$ and $P_2$ which exists, since $U_1\cap B$ is open and  connected and therefore connected by arc. If $\Gamma_C$ is the sub arc of $\Gamma$ linking $P_1$ and $P_2$, then the arcs $\Gamma_C$  and $\gamma_C$ have the same end points $P_1$ and $P_2$ and they are homotopic since $B$ is simply connected,
with the important point that $\gamma_C$ is now entirely contained in $U_1\cap B$.

Observe now that there are at most a countable number of connected components $C$ of $\Gamma\cap B.$ To see this we make a continuous parameterization $t\mapsto f(t)$ of the loop $\Gamma$ by $[0,1]$ such that $f(0)=f(1)\in U_1$. The set of $t$ such that $f(t)$ is in $B$ is an open subset of $(0,1)$ which is a finite or countable union of open interval, each generating a connected components $C$ of $\Gamma\cap B.$ Now we repeat the procedure for each $C$ and we land on a loop $\Gamma_1$ homotopic to
$\Gamma$ which is contained in $U_1$, thus homotopic to one point since $U_1$ is simply connected. Therefore $U_1\cup B$ is simply connected.

\vspace{4mm}\textsc{Continuation  of the proof of part 2.}

Consider the application $$m\mapsto h(m)=(h_1(m),\ldots,h_n(m))=\left(e^{\psi_1(m)},\ldots,e^{\psi_n(m)}\right)=e^{\psi (m)}$$ from $M_F$ to $(0,\infty)^n.$ Since $M_F$ is simply connected in $\R^n$ as the image of the convex set $\Theta(\mu)$ by the  diffeomorphism $\theta\mapsto k'(\theta)$ there exists an open subset $U_1 \subset \C^n$ simply connected such that $U_1\cap \R^n=M_F$
and a unique analytic application from $U_1$ to $\C^n$ whose  restriction to $ M_F$ is $h.$ For simplicity we still call $h$ this extension.

Observe that on $\Theta(\mu)\subset \R^n$ we have $$m_i=(k'_i(\theta)=\frac{\partial L}{\partial \theta_i}(\theta))/L(\theta)=e^{\theta_i}\frac{\partial f}{\partial z_i}(e^{\theta})/f(e^{\theta})$$ that we rewrite in $M_F$
$$m_i=e^{\psi_i(m)}\frac{\partial f}{\partial z_i}(z)/f(z)=h_i(m)/g_i(h(m)),$$ or in a more compact form the equation \eqref{LAGRANGE} but reduced to $M_F\subset (0,\infty)^n$

The Lagrange Theorem, or Proposition 6.1, implies the existence of an open ball $B\subset \C^n$ centered on 0 and the existence of an analytic function $\tilde{h}$ from $B$ to $U$ satisfying \eqref{LAGRANGE}. The ball  $B$ can be taken small enough such that $U_1\cap B$ is  connected and such that $0$ is the only zero of $\tilde{h}.$ Now we define $\mathcal{O}=U_1\cup B$. The fact of topology given above implies that $\mathcal{O}$  is simply connected. The two analytic functions $h$ defined on $U_1$ and  $\tilde{h}$ defined on $B$ coincide on the non empty open subset $B\cap U_1\cap M_F$  of $\R^n.$ Therefore they coincide on the connected subset $B\cap U_1$ and for simplicity we extend  the definition of $h$ to $\mathcal{O}$ by $h=\tilde{h}.$ Part 2 is proved.

\vspace{4mm}\textsc{Proof of part 3.}
$$P(m,F)(dx)=
\frac{e^{\psi (m)^T\, x}}
{L_{\mu}(\psi(m))}\mu(dx)
=\sum_{k\in \N^n}\mu_k\frac{(e^{\psi(m)})^k}
{L(\psi_{\mu}(m))}\delta_k(dx)=\sum_{k\in \N^n}\mu_k\frac{h(m)^k}
{f(h(m))}\delta_k(dx)$$
To conclude, we have
$$\psi_i(m)=\log h_i(m),\ \ \frac{\partial\psi_i}{\partial m_j}(m)=\frac{1}{h_i(m)}\times \frac{\partial h_i}{\partial m_j}(m)$$
In other terms
$$\psi'_i(m)=\mathrm{diag}(h_1(m)^{-1},\ldots,h_n(m)^{-1})\times h'(m),\ \ V(m)=(\psi'_i(m))^{-1}=(h'(m))^{-1}\mathrm{diag}
(h(m))$$ which is the desired result.

\vspace{4mm} The next proposition does not apply only on exponential families on $\N^n$ but is a technical result .

\vspace{4mm}\noindent\textbf{Proposition 6.4} Let $F(\mu)$ be an exponential family on $\R^n$ with  variance function $ V.$ Then there exists an open subset $U$ of $\C^n$ containing $M_F$ such that the function $\Phi(m)=k_{\mu}(\psi_{\mu}(m))$ defined on $M_F$ is extensible to $U$ as an  analytic function. This function is also denoted by $\Phi(m)$ and its restriction to $M_F$ satisfies $\Phi'(m)=(V(m))^{-1}m$ for all $m\in M_F.$

\vspace{4mm}\noindent\textbf{Proof.} Here again, for simplification we omit $\mu$   in $k_{\mu}$ and $\psi_{\mu}.$
Clearly we have $$(k(\psi(m)))'=\psi'(m)(k'(\psi(m)))=(V(m))^{-1}m.$$ Since $\psi$ is real-analytic on each point $m$ of $M_F$
there exists an open subset $\mathcal{O}_1\subset \C^n$ containing $M_F$ on which $\psi$ can be extended  and such that this extension $\psi$ is such that $\psi'$ has no zero on $\mathcal{O}_1.$ This grants that $\psi(\mathcal{O}_1)$ is also open. Similarly, we extend $\theta\mapsto k(\theta)$  to an open subset $U_1$ of $\C^n$ containing $\Theta(\mu).$ Since
$$ \psi (\mathcal{O}_1)\supset \psi(M_F)=\Theta(\mu)$$ we can claim that $U_2=U_1\cap \psi (\mathcal{O}_1)$ is a non empty open subset of $\C^n$. Finally we define $U= \psi^{-1}(U_2)$ and $\Phi(z)=k(\psi(z))$ on $U.$

\vspace{4mm}\noindent\textbf{Proof of Theorem 6.2.}  We prove first the 'if' part as well as formulas \eqref{FUNCTIONG} and \eqref{LAGRANGEMU}. Therefore we assume that properties 1) and 2) are true.

From Proposition 6.4 there exists $U\subset\C^n$ such that $U\supset M_F$ and there exists $\Phi$ analytic on $U$
such that for $m$ in $M_F$, we have $\Phi'(m)=(V(m))^{-1}m.$ Let $W$ be the connected component of $\mathcal{O}\cap U$ containing $M_F$. Such a $W$ does exist since $M_F$ is connected. Therefore the restriction of $\Phi'$ to $W$ coincides with $\phi'$, since they are analytic and coincide on $M_F$ which is an open set of $\R^n.$ Now if $z\in \mathcal{O}$ and
$$\phi'(z)=(\phi'_1(z),\ldots,\phi'_n(z))^T,$$  then for each $i,j$ in $\{1,\ldots,n\}$ the map
$$z\mapsto \frac{\partial \phi_i}{\partial z_j}-\frac{\partial \phi_j}{\partial z_i}$$ is analytic on $\mathcal{O}$ and is zero on $W$, and  as a consequence of connectivity, also on $\mathcal{O}$. The local condition of gradient is verified for $\phi'.$  Since $\mathcal{O}$ is simply connected, $\phi'$ is the gradient of a function $\phi$ defined on $\mathcal{O}.$

Let us now show that there exists an unique analytic function $K$ defined in a neighborhood $\mathcal{V}\subset \C^n$ of $0$ and valued in $\C^n$ such that $K'(0)=0$ and such that for all $z$ in $\mathcal{V}$ one has
$$K'(z)(z)=\mathbf{1}-\phi'(z).$$

First there  exists $\mathcal{V}\subset \C^n$ such that for $z\in \mathcal{V}$ we have
$$\phi(z)=\sum_{k\in \N^n}b_kz^k$$ for some $(b_k)_{k\in \N^n}.$ The function $K$ that we are looking  for should satisfy for all $i=1,\ldots,n$
\begin{equation}\label{THEBIGK}\sum_{j=1}^nz_j\frac{\partial K_i}{\partial z_j}=1-\frac{\partial \phi}{\partial z_i}\end{equation}
We write $K_i$ as a sum of a series:
$$K_i(z)=\sum_{k\in \N^n}a_k^iz^k$$
With this notation \eqref{THEBIGK} is equivalent to
$$ \sum_{k\in \N^n}|k|a_k^iz^k=1-\frac{\partial \phi}{\partial z_i}.$$ Therefore $a_k^i$ is computable in terms of
$(b_k)_{k\in \N^n}$  if $|k|\neq 0.$ For computing $a_0^i$ we use the fact that $K_i(0)=\phi_i(0)-1=0,$ implying
$a_0^i=0.$ Hence the function $K$ is completely known and unique.

We are now in position to construct the function $G$ appearing  on \eqref{FUNCTIONG}. It is $G:\mathcal{V}\rightarrow
\C^n$
$$G(z)=(e^{K_1(z)},\ldots, e^{K_n(z)})=e^{K(z)}.$$
It satistisfies
$$G'(z)=\mathrm{diag}(e^{K(z)})K'(z)=\mathrm{diag}(G(z))K'(z)$$
leading to

$$\mathrm{diag}\left(\frac{1}{G_1(z)},\ldots,\frac{1}{G_n(z)}\right)\times G'(z)=K'(z)$$
Applying the last line to $z$ we get

$$\mathrm{diag}\left(\frac{1}{G_1(z)},\ldots,\frac{1}{G_n(z)}\right)\times G'(z)z=K'(z)z$$
and  using \eqref{THEBIGK} we get \eqref{FUNCTIONG}

$$\mathrm{diag}\left(\frac{1}{G_1(z)},\ldots,\frac{1}{G_n(z)}\right)\times G'(z)z=\mathbf{1}-
\phi'(z).$$ and, since $K(0)=0$ we get $G(0)=\mathbf{1}.$

In order to show \eqref{LAGRANGEMU} we now consider for $m\in M_F$ an arbitrary measure $\nu$ generating $F$
(Therefore $F=F(\mu)=F(\nu)),$ as well the following function $H.$
$$H(m)=\mathrm{diag}\, (m) e^{-\psi_{\nu}(m)}.$$ We are going to prove that there exists positive constants $c=(c_1,\ldots,c_n)$ such that \begin{equation}\label{FUNCTIONH}H(m)=\diag(c)G(m).\end{equation}
With obvious notations we write
$$\diag (H(m)^{-1})=\mathrm{diag}\, (m^{-1})\diag\left( e^{\psi_{\nu}(m)}\right).$$  Compute $H'(m)$ as follows
\begin{eqnarray*}
H'(m)&=&\diag (e^{-\psi_{\nu}(m)})-\diag (m)]\diag(e^{-\psi_{\nu}(m)})\psi'_{\nu}(m)\\
\diag (H(m)^{-1})H'(m)&=&\mathrm{diag}\, (m^{-1})\diag\left( e^{\psi_{\nu}(m)}\right)\left(\diag (e^{-\psi_{\nu}(m)})-\diag (m)\diag((e^{-\psi_{\nu}(m)})\psi'_{\nu}(m)\right)\\
&=&\mathrm{diag}\, (m^{-1})\left(I_n-\diag (m)\psi'_{\nu}(m)\right)=\mathrm{\diag}\, (m^{-1})\left(I_n-\diag (m)V(m)^{-1}\right)\\
\diag (H(m)^{-1})H'(m)(m)&=&\mathrm{diag}\, (m^{-1})\left(I_n-\diag (m)V(m)^{-1}\right)m\\&=&\mathbf{1}-V(m)^{-1}m=\mathbf{1}-\phi'(m)
\end{eqnarray*}
Therefore $H$ satisfies also \eqref{FUNCTIONG}. For convenience write $H= e^L$ (since $M_F\subset (0,\infty)^n$ it makes sense). Recall that $G=e^K$. The fact that both $G $ and $H$ satisfy \eqref{FUNCTIONG} imply that
$$\diag (m )K'(m)=\diag (m)L'(m),\ \diag (m )(K'(m)-L'(m))=0,$$ or for all $i=1,\ldots,n$
$$\sum_{j=1}^nm_j\frac{\partial(K_i(m)-L_i(m))}{\partial m_j}=0.$$ This shows that the function $K_i(m)-L_i(m)$ is homogeneous of degree $0$. Therfore this function is constant along rays from the origin. In other words, it has the form $f(m/\|m\|)$ which cannot be analytic on $0$ unless $f$ is a constant on the whole $\R^n.$. This proves \eqref{FUNCTIONH}.

Let us apply the Lagrange formula of Proposition 6.1 to $g=G$. Observe that $D(G)$ as defined by \eqref{JACOBIAN2} is not zero in a neighborhood of  $0$ in $\C^n$. Actually \eqref{JACOBIAN1} says
$$M(G)=I_n-\diag (\frac{z_1}{G_1},\ldots,\frac{z_n}{G_n})G'(z)=I_n-\diag(\frac{z}{G(z)})G'(z)$$
Since $G=e^K$ and since $K$ is analytic in a neighborhood of $0$ then $\lim_{z\to 0} M(G(z))=I_n$ which implies
the existence of a neighorhood of zero such that $D(G)\neq 0$ on it.

From Proposition 6.1 there exists a neighborhood $ \mathcal{W}$ of 0 and a function $h$ from $ \mathcal{W}$ to
$ \mathcal{V}$ such that
$$h(w)=\diag (w)G(h(w))$$
We also use for $g_0$ the function $g_0=e^{\phi}.$

From Proposition 6.1 again, and with $(\mu_k)_{k\in \N^n}$ defined by \eqref{LAGRANGEMU} we can write for $w\in \mathcal{W}$
$$
e^{\phi(h(w))}=\sum_{k\in \N^n}\mu_kw^k.
$$
Also for $ m=h(w)\in M_F:$\begin{equation}\label{LAGRANGEMU3}
e^{\phi(m)}=\sum_{k\in \N^n}\mu_kw^k.
\end{equation}
We now try to get an expression of $w$ in terms of $m$ in formula \eqref{LAGRANGEMU3}.
$$m=h(w)=\diag (w) G(m)=\diag (w) \diag (c) H(m)=\diag (w) \diag (c) \diag( m))e^{-\psi_{\nu}(m)}$$
Watching this formula for each component we obtain
$$m_i=w_ic_im_ie^{-(\psi_{\nu}(m))_i}\ \ \Rightarrow w_i=c_i^{-1}e^{(\psi_{\nu}(m))_i}$$Pligging in \eqref{LAGRANGEMU3} we get:
\begin{equation}\label{LAGRANGEMU4}
e^{\phi(m)}=\sum_{k\in \N^n}\mu_k(c^{-1})^ke^{\langle k,\psi_{\nu}(m)\rangle}.\end{equation}
We come back to $\phi'(m)=V(m)^{-1}(m)=\psi'_{\nu}(m)m$ giving $\phi(m)=k_{\nu}(\psi_{\nu}(m))+b$ where $b$ is a constant. This comes from the fact that $\nu$ is an arbitrary generating measure of $F$.

As a result in terms of Laplace transforms we have
$$e^{\phi(m)}= e^bL_{\nu}(\psi_{\nu}(m))=L_{e^b\nu}(\psi_{\nu}(m))$$
Comparing with \eqref{LAGRANGEMU4} we get that
$$\nu=\sum_{k\in \N^n}\mu_ke^{-b}(c^{-1})^k\delta_k=\sum_{k\in \N^n}\nu_k\delta_k.$$ Therefore $\mu_k =e^{<k,-\log c>+b}\nu_k$
and this says that $\mu=\sum_{k\in N^n}\mu_k\delta_k$ is also a basis of $F$. This ends the 'if' part of Theorem 6.2.

We now prove the 'only if' part. We assume that $F(\mu)$ is concentrated on $\N^n$ with $\mu_0,\mu_{e_1},\ldots,\mu_{e_n}$ positive and we want to show 1) and 2). Actually this is a direct consequence of Proposition 6.3. Part1 of Theorem 6.2 is given in Part 1 of Proposition 6.3. For Part 2 the existence of $\mathcal{O}$ is in Part2 of Proposition 6.3.

The open set $\mathcal{O}$ mentioned in Proposition 6.3 contains $M_F$ and $0$. Therefore $\mathcal{O}\cap U$ is not empty. Denote by $W$ the connected component of $\mathcal{O}\cap U$ containing $M_F$ and $0$ (which does exist since $M_F$ is connected and $0\in\overline{M_F}$). The restriction of $\phi$ to $W$ coincides with $\Phi'$ since they coincide on $M_F$, which is an open subset of $\R^n.$

Finally the function $\phi'(m)=\diag(\frac{1}{h(m)})h'(m)m$ is $(V(m))^{-1}m$ from Part 3. Last thing to prove is that $\phi'(0)=\mathbf{1}.$

To do this, we apply Proposition 6.1 to the function $g$ defined in Proposition 6.3
part 2, and to the function $g_0$ as the projection map $(z_1,\ldots,z_n)\mapsto z_i.$ We get

$$h_i(m)=\sum_{k\in \N^n}m^n[z^n](z_ig{z}^nD(g_1,\ldots,g_n)(z))$$
In this last expression, we use the fact that $g_i(z)=f(z)/\frac{\partial f}{\partial z_i}$ and
$$f(z)=\mu_0+\mu_{e_1} z_1+\ldots+\mu_{e_n} z_n+\sum_{k\in \N^n,|k|\geq 2}\mu_k z^k$$
for computing $$\frac{\partial h_i}{\partial m_j}(0)=[z_j](z_ig(z)^nD(g_1,\ldots,g_n)(z))$$
Observe that this number is $0$ if $i\neq j$ and that
$$\frac{\partial h_i}{\partial m_i}(0)=[z_i](z_ig_i(z)D(g_1,\ldots,g_n)(z))=g_i(0)=\frac{\mu_0}{\mu_{e_i}}$$
since $g_i(z)=f(z)/\frac{\partial f (z)}{\partial z_i}$ and that $D(g)(0)=1.$
Coming back to
$$\phi'_i(m)=\frac{1}{h_i(m)}\sum_{j=1}^n   \frac{\partial h_i(m)}{\partial m_j}\times m_j$$
we do $m=m_ie_i$ in the above expression for getting
$$\phi'_i(m_ie_i)=\frac{1}{h_i(m_ie_i)}\frac{\partial h_i(m_ie_i)}{\partial m_i}\times m_i$$
Since we have seen, using the hypothesis on $ \mu$, that $\frac{\partial h_i}{\partial m_i}(0)=\frac{\mu_0}{\mu_{e_i}}$ exists and is not zero
we can claim by letting $m_i\to 0$ that $\phi'_i(0)=1.$

The proof of  Theorem 6.2 is finished.

\vspace{4mm}

\vspace{4mm}\noindent\textbf{Following the proof  of Theorem 6.2 through an example.}
Let
$$V(m)=\mathrm{diag}(m)+m\otimes m$$
We decide to follow step after step the proof of Theorem 6.2 for $(\mu_k)_{k\in \N^n}.$ For this particular case, there are other ways to compute $\mu$ (see Letac \cite{Letac}). One arrives easily to
$$\Phi'(m)=\frac{1}{1+s}\mathbf{1}$$ with  $s=m_1+\cdots +m_n.$
Now we have to compute  $K$ appearing in the proof of Theorem 6.2 with the only  information $$K'(m)m=\mathbf{1}-\Phi'(m)$$ translated in $n$ equations
 $i=1,\ldots ,n$ as follows
$$\sum_{k=1}^nm_k\frac{\partial}{\partial m_k}K_i(m)=1-\Phi_i'(m).$$ Even with the condition $\Phi'(0)=\mathbf{1}$ or $K(0)=0$, this does not seem enough for computing $K$ since the homogeneous equation
$$\sum_{k=1}^dm_k\frac{\partial}{\partial m_k}K_i(m)=0$$
has many solutions. However, these solutions have the form $K_i(m)=f_i(m/\|m\|)$ which cannot be analytic on $ 0$
if $f_i$ is not a constant function. Fortunately, analyticity of $K$ allows us to write
$$\sum_{k=1}^nm_k\frac{\partial}{\partial m_k}K_i(m)=\sum_k(k_1+\cdots+k_n)a_k^im^k=\frac{s}{1+s}$$ which gives us in theory to compute the $a^i_k$  -except $a^i_0$!
However $K(0)=0$ gives finally

$$K_i(m)=\sum_{k\neq 0}\frac{(-1)^{|k|-1}}{|k|}\left( \begin{array}{c} |k|\\k\end{array}\right)m^k=\log (1+s)$$ hence $$G(m)=(1+s)\mathbf{1}.$$ We obtain $\mu_k$ for $k\neq 0:$
\begin{eqnarray*}\mu_k&=&[m^k]e^{\Phi(m)}G(m)^k \det \left(I_d-(\frac{m_j}{G_i(m)}\frac{\partial G_i(m)}{\partial m_j})\right)\\&=&[m^k]((1+s)(1+s)^k \mathbf{1}\det \left(I_n-\frac{m\otimes \mathbf{1}}{1+s}\right)\\&=&[m^k]((1+s)^{[k]}=\left(\begin{array}{c}|k|\\k\end{array}\right)\end{eqnarray*}
 with the usual formula $\det (I_d-a \otimes b)=1-\langle a,b\rangle.$
As a consequence
$$ \sum_{k\in \N^n} \mu_kz^k=\sum_{j=0}^{\infty}(z_1+\cdots+z_n)^j=\frac{1}{1-z_1-\cdots-z_n}.$$

\section{The action of the group $H$ on absolutely continuous or discrete exponential families}
 Recall that the subgroup $H$ of $G$ has been defined in \eqref{GROUPEH} and is made with the matrices   of the form $g=\left[\begin{array}{cc}I_n&0\\c^T&\lambda\end{array}\right]$ where $\lambda>0$ and $c\in \R^n$. In this section we study the action of these $T_g$'s on the exponential families through their generating measures and the Laplace transforms. The important point is that $H$ preserves the absolutely continuity (Proposition 7.2)  or the set of measures concentrated on $\Z^d$ (Proposition 7.4). The basis of the proof is a stricking result due to Philippe Rouqu\`es \cite{Rouqu} in his unpublished Master Thesis in Toulouse, where he  gives light  to and extends the results of Terry Speed \cite{Speed}. The bulk of his reasoning is in Proposition 7.3. Our Proposition 7.1 is an adaptation of it to the absolutely continous case. The proofs of Propositions 7.1 and 7.3 are quite similar, but we gave up the idea of gathering them in one statement for clarity.

\vspace{4mm}\noindent\textbf{Proposition 7.1.}  Let $(\mu_{\lambda}(dx))_{\lambda>0}$ a convolution semigroup on $\R^n$, such that $ 0\in \Theta(\mu_{\lambda }).$ We assume that  $\mu_{\lambda}(dx)$ has density $f(\lambda,x)$ and we consider the density on $\R^n$ defined  for $c$ in $\R^n$ by
\begin{equation}\label{ROUQUES3}
p(\lambda,c,x)=\frac{\lambda}{\lambda+\<c,x\>}f(\lambda+\<c,x\>,x)1_{(0,\infty)}(\lambda+\<c,x\>).
\end{equation}
Then if $c\neq 0$ \begin{equation}\label{ROUQUES4}p(\lambda,c,x)1_{(0,\infty)}(\<c,x)\>)=\left(\int_{\R^n}p(\lambda,0,y)p(\<c,y\>,c,x-y)dy\right).\end{equation}

\vspace{4mm}\noindent\textbf{Proof.} It is useful to reformulate \eqref{ROUQUES4} by observing that
$$p(\<c,y\>,c,x-y)=\frac{\<c,y\>}{\<c,x\>}f (\<c,x\>,x-y )1_{(0,\infty)}(\<c,x\>)$$ so we have to prove that
\begin{equation}\label{ROUQUES5}p(\lambda,c,x)1_{(0,\infty)}(\<c,x\>)=\frac{1}{\<c,x\>}\left(\int_{\R^n}f(\lambda, y)\<c,y\>f(\<c,x\>,x-y )dy\right)1_{(0,\infty)}(\<c,x\>)\end{equation}

 Denote  $L_{\lambda}(\theta)=(L_1(\theta))^{\lambda}=\int_{\R^n}e^{\<\theta,y\>}f(\lambda,y)dy.$ Let us fix a number $A>0$, consider the differential operator $D=\sum_{i=1}^n c_i\frac{\partial}{\partial \theta_i}$ and observe that
\begin{equation}\label{ROUQUES6}L_{A\lambda}\times  D\,  L_{\lambda}=\frac{1}{A+1}D\, L_{(A+1)\lambda}.\end{equation}
This comes from  $$D L_{\lambda}=D (L_1^{\lambda})=\lambda L_1^{\lambda-1}D (L_1)\Rightarrow
L_{A\lambda}D L_{\lambda}=\lambda L_1^{(A+1)\lambda-1}D (L_1)=\frac{\lambda}{(A+1)\lambda}D(L^{(A+1)\lambda})$$
The right  hand side of \eqref{ROUQUES6} is $$\frac{1}{A+1}\int_{\R^n}e^{\<\theta,y\>}\<c,y\>f((A+1)\lambda,y)dy$$

The left hand side of \eqref{ROUQUES6} is the product of the  Laplace transforms of the two  respective densities
$f(A\lambda,y)$ and $\<c,y\>f(\lambda,y)$. Since $ 0\in \Theta(\mu_{\lambda })$ this implies that $y
\mapsto \<c,y\>f(\lambda,y)$ is integrable. Therefore

$$\int_{\R^n}\<c,y\>e^{\<\theta,y\>}f(\lambda,y)dy\times\int_{\R^n}e^{\<\theta,y\>}f(A\lambda,y)dy
=\frac{1}{A+1}\int_{\R^n}e^{\<\theta,y\>}\<c,y\>f((A+1)\lambda,y)dy.$$

As a result we can claim

$$\frac{1}{A+1}\<c,x\>f((A+1)\lambda,x)=\int_{\R^n}\<c,y\>f(\lambda,y)f(A\lambda,x-y)dy$$

Now in this formula we replace $A$ by $\<c,x\>/\lambda$, we multiply both sides by $1_{(0,\infty)}(\<c,x)\>) $ and we get exactly \eqref{ROUQUES4} since $1_{(0,\infty)}(
\lambda+\<c,x)\>) 1_{(0,\infty)}(\<c,x)\>) =1_{(0,\infty)}(\<c,x)\>) .$

\vspace{4mm}\noindent\textbf{Example.} Let us perform the calculations of the proof of Proposition 7.1  for the particular case $n=1$ and for $$\mu_{\lambda}(dx)=N(0,\lambda)(dx)=f(\lambda,x)dx=e^{-\frac{x^2}{2\lambda}}\frac{dx}{\sqrt{2\pi\lambda}}.$$
Then $$ f(\lambda+cx,x)=e^{-\frac{x^2}{2\lambda+cx}}\frac{1}{\sqrt{2\pi(\lambda+cx)}},\ \  \ p(\lambda,c,x)=\frac{\lambda}{\lambda+cx}e^{-\frac{x^2}{2(\lambda+cx)}}\frac{1}{\sqrt{2\pi(\lambda+cx)}}1_{(0,\infty)}(\lambda+cx).$$  Thus we get
$$\ p(cy,c,x-y) = \frac{cy}{cx}\frac{1}{\sqrt{2\pi cx}}e^{-\frac{(x-y)^2}{2cx}}1_{(0,\infty)}(cx)     $$
and equality \eqref{ROUQUES4} becomes
\begin{equation}\label{ROUQUES7}p(\lambda,c,x)1_{(0,\infty)}(cx)=\left(\int_{\R}ye^{-\frac{y^2}{2\lambda}}\frac{1}{\sqrt{2\pi\lambda}}\times \frac{1}{\sqrt{2\pi cx}}e^{-\frac{(x-y)^2}{2cx}}dy\right)\frac{1}{x}1_{(0,\infty)}(cx)\end{equation}
Now consider  the right hand side of \eqref{ROUQUES7}. We undertake its computation.  Let $x$ be a fixed constant in \eqref{ROUQUES7} and consider the new function of $z:$ $$g(z)=\int_{\R}ye^{-\frac{y^2}{2\lambda}}\frac{1}{\sqrt{2\pi\lambda}}\times \frac{1}{\sqrt{2\pi cx}}e^{-\frac{(z-y)^2}{2cx}}dy$$ which is nothing but the density of the convolution of $y\mapsto yf(\lambda,y)$ and of $y\mapsto f(cx,y).$

The Laplace transform of $y\mapsto yf(\lambda,y)$ is
\begin{equation}\label{ROUQUES8}
\int_{\R}ye^{-\frac{y^2+\theta y}{2\lambda}}\frac{1}{\sqrt{2\pi\lambda}}=\frac{d}{d\theta}\int_{\R}e^{-\frac{y^2+\theta y}{2\lambda}}\frac{1}{\sqrt{2\pi\lambda}}=\frac{d}{d\theta}e^{\frac{\lambda\theta^2}{\lambda}}=\lambda \theta e^{\frac{\lambda\theta^2}{2}}\end{equation} and of course the Laplace transform of
$f(cx,y)=e^{\frac{cx\theta^2}{2}}.$  The product of these two Laplace transforms, i.e. the Laplace transform of $g(z)$ is therefore $$\lambda \theta e^{\frac{(\lambda+cx)\theta^2}{2}}=\frac{\lambda}{\lambda+cx}\times \left((\lambda+cx)\theta e^{\frac{(\lambda+cx)\theta^2}{2}}\right)$$ From \eqref{ROUQUES8}, replacing $\lambda$ by $\lambda+cx$ we see that $$g(z)=\frac{\lambda}{\lambda+cx}\times\left(e^{-\frac{z^2}{2(\lambda+cx)}}\frac{1}{\sqrt{2\pi(\lambda+cx)}}\right)=\frac{\lambda }{\sqrt{2\pi}(\lambda+cx)^{3/2}}e^{-\frac{z^2}{2(\lambda+cx)}}$$
Finally we get the explicit expression of \eqref{ROUQUES7} by replacing $z$ by $x$
$$p(\lambda,c,x)1_{(0,\infty)}(cx)=\frac{\lambda  }{\sqrt{2\pi}(\lambda+cx)^{3/2}}e^{-\frac{x^2}{2(\lambda+cx)}}1_{(0,\infty)}(cx),$$ which fits with the definition of $p(\lambda,c,x).$

\vspace{4mm}\noindent\textbf{Proposition 7.2.} With the notations and hypotheses of Proposition 7.1, consider
$$\mu_{\lambda,c}(dx)=p(\lambda,c,x)1_{(0,\infty)}(\<c,x)\>) dx.$$ Then
\begin{enumerate}\item $\mu_{\lambda,c}$ is in $\mathcal{M}(\R^n)$
\item Its cumulant transform satisfies in a neighborhood of zero $$k_{\mu_{\lambda,c}}(\theta)
=k_{\mu_{\lambda}}(\theta+ck_{\mu_{1,c}}(\theta))$$

\item In particular if $g=\left[\begin{array}{cc}I_n&0\\c^T&\lambda\end{array}\right]$ is in the group $H$ defined in \eqref{GROUPEH}, we have $$T_g(F(\mu_1))=F(\mu_{\lambda,c})$$

\item $\mu_{\lambda,c}$ is infinitely divisible, with $\lambda$ as the Jorgensen parameter. Also, if $\mu_1$  is a probability then $\mu_{1,c}$ is also a probability.
\end{enumerate}

\vspace{4mm}\noindent\textbf{Proof.}

\begin{eqnarray} \nonumber
L_{\mu_{\lambda,c}}(\theta)&=&\int_{\R^n}e^{\<\theta,x\>}p(\lambda,c,x)1_{(0,\infty)}(\<c,x)\>) dx\\
&=&\label{FROM}\int_{\R^{2n}}e^{\<\theta,x\>}p(\lambda,0,y)p(\<c,y\>,c,x-y)dxdy\\
&=&\label{FROM1}\int_{\R^{n}}e^{\<\theta,y\>}p(\lambda,0,y)\left(\int_{\R^{n}}e^{\<\theta,z\>}p(\<c,y\>,c,z)dz\right)dy\\
&=&\nonumber\int_{\R^{n}}e^{\<\theta,y\>}p(\lambda,0,y)\left(L_{\mu_{1,c}}(\theta)\right)^{\<c,y\>}dy\\
&=&\nonumber\int_{\R^{n}}e^{\<\theta,y\>}p(\lambda,0,y)e^{k_{\mu_{1,c}}(\theta)\<c,y\>}dy\\
&=&\nonumber\int_{\R^{n}}e^{\<\theta+ck_{\mu_{1,c}},y\>}p(\lambda,0,y)dy=L_{\mu_{\lambda}}(\theta+ck_{\mu_{1,c}}(\theta))
\end{eqnarray}
Line \eqref{FROM} comes from \eqref{ROUQUES4} and line \eqref{FROM1} comes from the change of variable $z=x-y.$ This proves both Part 1 and Part 2.

For proving Part 3 we just apply Proposition 3.3 by replacing $(\mu,\nu,
\lambda)$ by $(\mu_{\lambda,c}, \mu_1,\theta).$
From Part 2 we have $k_{\mu_{\lambda,c}}(\theta)
=k_{\mu_{\lambda}}(\theta+ck_{\mu_{\lambda,c}}(\theta)).$
Denote $\alpha(\theta)=\theta+ck_{\mu_{\lambda,c}}(\theta)$ and therefore

$$\left[\begin{array}{c}\theta\\-k_{\mu_{\lambda,c}}(\theta)\end{array}\right]=\left[\begin{array}{cc}I_n&c\\0&\lambda\end{array}\right]\left[\begin{array}{c}\alpha(\theta)\\-k_{\mu_1}(\alpha(\theta))\end{array}\right]
$$
Here the function $\theta(\lambda)$ appearing in the statement of Proposition 3.3 is replaced by $\alpha(\theta).$
Infinite divisibility mentioned in Part 4 is obvious from Part 2. Also
$$\left[\begin{array}{c}\alpha\\-k_{\mu_{1}}(\alpha)\end{array}\right]=\left[\begin{array}{cc}I_n&-c\\0&\lambda\end{array}\right]\left[\begin{array}{c}\theta\\-k_{\mu_{\lambda,c}}(\theta)\end{array}\right]
$$ As a consequence $k_{\mu_{1}}(\alpha)
=k_{\mu_{1,c}}(\alpha-ck_{\mu_{1}}(\alpha)).$ If $\mu_1$ is a probability then $k_{\mu_{1}}(0)=0$ and clearly
$k_{\mu_{1,c}}(0)=0.$

\vspace{4mm}\noindent\textbf{Example: $\mu_{\lambda}$ generates a Gamma family.}
Let us apply Proposition 7.1 to $n=1$ and to $$f(\lambda,x)=\frac{x^{\lambda-1}}{\Gamma(\lambda)}1_{(0,\infty)}(x).$$ Note that $\Theta(\mu)=(-\infty,0).$ Its Laplace transform is $L_{\mu_{\lambda}}(\theta)=(-\theta)^{-\lambda}.$ We intend to compute the density of $\mu_{\lambda,c}(dx)$ and to say something about its Laplace transform.  For simplicity we assume $c>0.$
With the notations of Propositions 7.1 and 7.2 the density of $\mu_{\lambda,c}(dx)$ is
$$p(\lambda,c,\lambda)1_{(0,\infty)}(cx)=\frac{\lambda}{\lambda+cx}\frac{x^{\lambda+cx}}{\Gamma(\lambda+cx)}1_{(0,\infty)}(x).$$
Now from Part 2 of Proposition 7.2, we have
$$k_{\mu_{\lambda,c}}(\theta)=-\lambda\log (-\theta+k_{\mu_{1,c}}(\theta)). $$
For simplicity denote $z(\theta)=k_{\mu_{1,c}}(\theta)$. It satisfies
$$L_{\mu_{\lambda,c}}(\theta)=e^{\lambda z(\theta)}=\frac{1}{(-\theta+cz(\theta))^{\lambda}}.$$
The functional equation $e^{-z(\theta)}=cz(\theta)-\theta$ for computing $z(\theta)$ can be treated by the Lagrange formula. To conclude this example, since the variance function of $F(\mu_{\lambda})$ is $m^2/\lambda,$ according to the formula \eqref{TG} we have, since $h_g(m)=m/(\lambda+cm):$
$$V_{F(\mu_{\lambda,c})}(m)=T_g(V_{F(\mu_{\lambda})})(m)=\frac{1}{\lambda^3}m^2(\lambda+cm)$$
namely the variance function of a Kendall-Ressel family.

\vspace{4mm}\noindent\textbf{Proposition 7.3.} Let $(\mu_{\lambda}(dx))_{\lambda>0}$
a convolution semigroup on $\Z^n$. Let $c\in \R^n$. Denote
$$\mu_{\lambda}(dx)=\sum_{k\in \Z^n}f(\lambda,k)\delta_k(dx)$$
Define also $$p_k(\lambda,c)=\frac{1}{\lambda+\<c,k\>}f(\lambda+\<c,k\>,k)1_{(0,\infty)}(\lambda+\<c,k\>).$$
Then \begin{equation}\label{ROUQUES5}p_k(\lambda,c)1_{[0,\infty)}(\<c,k\>)=\sum_{k'\in \Z^n}p_{k'}(\lambda,0)p_{k-k'}(\<c,k'\>,c)\end{equation}

\vspace{4mm}\noindent\textbf{Proof.} It is useful to reformulate \eqref{ROUQUES5} by observing that
$$p_{k-k'}(\<c,k'\>,c)=\frac{\<c,k'\>}{\<c,k\>}f (\<c,k\>,k-k' )1_{(0,\infty)}(\<c,k\>)$$ so we have to prove that
\begin{equation}\label{ROUQUES51}p_k(\lambda,c)1_{(0,\infty)}(\<c,k\>)=\frac{1}{\<c,k\>}\left(\sum_{k'\in
\Z^n}\<c,k'\>f(\lambda, k')f(\<c,k\>,k-k')\right)1_{(0,\infty)}(\<c,k\>))\end{equation}

 Denote  $L_{\lambda}(\theta)=(L_1(\theta))^{\lambda}=\sum_{k'\in
\Z^n}e^{\<\theta,k'\>}f(\lambda,k').$ Let us fix a number $A>0$, consider the differential operator $D=\sum_{i=1}^n c_i\frac{\partial}{\partial \theta_i}$ and observe that
\begin{equation}\label{ROUQUES61}L_{A\lambda}\times  D\,  L_{\lambda}=\frac{1}{A+1}D\, L_{(A+1)\lambda}.\end{equation}
This comes from  $$D L_{\lambda}=D (L_1^{\lambda})=\lambda L_1^{\lambda-1}D (L_1)\Rightarrow
L_{A\lambda}D L_{\lambda}=\lambda L_1^{(A+1)\lambda-1}D (L_1)=\frac{\lambda}{(A+1)\lambda}D(L^{(A+1)\lambda})$$
The right  hand side of \eqref{ROUQUES61} is $$\frac{1}{A+1}\sum_{k'\in
\Z^n}e^{\<\theta,k'\>}\<c,k'\>f((A+1)\lambda,k')$$

The left hand side of \eqref{ROUQUES61} is the product of the  Laplace transforms of the two  respective measures on $\Z^n$
$$\sum_{k'\in \Z^n}f(A\lambda,k')\delta_{k'},\ \ \  \sum_{k\in \Z^n}\<c,k\>f(\lambda,k)\delta_k.$$ Since $ 0$ is in $\Theta(\mu_{\lambda })$ this implies that $k
\mapsto \<c,k\>f(\lambda,k)$ is summable. Therefore

$$\sum_{k'\in
\Z^n}\<c,k'\>e^{\<\theta,k'\>}f(\lambda,k')\times\sum_{k'\in
\Z^n}e^{\<\theta,k'\>}f(A\lambda,k')
=\frac{1}{A+1}\sum_{k'\in
\Z^n}e^{\<\theta,k'\>}\<c,k'\>f((A+1)\lambda,k')$$

As a result we can claim

$$\frac{1}{A+1}\<c,k\>f((A+1)\lambda,k)=\sum_{k'\in
\Z^n}\<c,k'\>f(\lambda,k')f(A\lambda,k-k')$$

Now in this formula, we replace $A$ by $\<c,k\>/\lambda$ , we multiply both sides by $1_{(0,\infty)}(\<c,k)\>) $ and we get exactly  \eqref{ROUQUES5} since $1_{(0,\infty)}(
\lambda+\<c,k)\>) 1_{[0,\infty)}(\<c,k)\>) =1_{(
[0,\infty)}(\<c,k)\>) .$

\vspace{4mm}\noindent\textbf{Comment.} For simplicity, in Proposition 7.3  we have dealt with infinitely divisible measures. We could have  written a slight generalisation by dealing with  a Jorgensen set different from $(0,\infty)$ for instance with
$\Lambda(\mu)$ equal to the set of positive integers. This is indeed in the Bernoulli case $\mu=\delta_0+\delta_{e_1}+\cdots+\mu_{e_n}.$ This implies that to have $c\in \R^n$ such that $\<c,k\>$ in $\N$. In the Bernoulli case such a modified Proposition 7.3
is exactly formulas [2] and [2'] of \cite{Speed}.

\vspace{4mm}\noindent\textbf{Proposition 7.4.} With the notations and hypotheses of Proposition 7.3, consider

$$\mu_{\lambda,c}(dx)=\sum_{k\in \Z^n}p_k(\lambda,c)1_{(0,\infty)}(\<c,k\>) \delta_k(dx).$$ Then
\begin{enumerate}\item $\mu_{\lambda,c}$ is in $\mathcal{M}(\R^n)$
\item Its cumulant transform satisfies in a neighborhood of zero $$k_{\mu_{\lambda,c}}(\theta)
=k_{\mu_{\lambda}}(\theta+ck_{\mu_{1,c}}(\theta))$$

\item In particular if $g=\left[\begin{array}{cc}I_n&0\\c^T&\lambda\end{array}\right]$ is in the group $H$ defined in \eqref{GROUPEH}, we have $$T_g(F(\mu_1))=F(\mu_{\lambda,c})$$

\item $\mu_{\lambda,c}$ is infinitely divisible, with $\lambda$ as the Jorgensen parameter. Also, if $\mu_1$  is a probability then $\mu_{1,c}$ is also a probability.
\end{enumerate}

\vspace{4mm}\noindent\textbf{Proof.}

\begin{eqnarray} \nonumber
L_{\mu_{\lambda,c}}(\theta)&=&\sum_{k\in \Z^n}e^{\<\theta,k\>}p_k(\lambda,c)1_{(0,\infty)}(\<c,k)\>) \\
&=&\label{GROM}\sum_{(k,k')\in \Z^{2n}}e^{\<\theta,k\>}p_{k'}(\lambda,0)p_{k-k'}(\<c,y\>,c)\\
&=&\label{GROM1}\sum_{k'\in \Z^n}e^{\<\theta,k'\>}p_{k'}(\lambda,0)\left(\sum_{m\in \Z^n}e^{\<\theta,m\>}p_m(\<c,k'\>,c)\right)\\
&=&\nonumber \sum_{k'\in \Z^n}e^{\<\theta,k'\>}p_{k'}(\lambda,0)\left(L_{\mu_{1,c}}(\theta)\right)^{\<c,k'\>}\\
&=&\nonumber \sum_{k'\in \Z^n}e^{\<\theta,k'\>}p_{k'}(\lambda,0)e^{k_{\mu_{1,c}}(\theta)\<c,k'\>}\\
&=&\nonumber \sum_{k'\in \Z^n}e^{\<\theta+ck_{\mu_{1,c}},k'\>}p_{k'}(\lambda,0)=L_{\mu_{\lambda}}(\theta+ck_{\mu_{1,c}}(\theta))
\end{eqnarray}
Line \eqref{GROM} comes from \eqref{ROUQUES5} and line \eqref{GROM1} comes from the change of variable $m=k-k'.$ This proves both Part 1 and Part 2.

For proving Part 3 we just apply Proposition 3.3 by replacing $(\mu,\nu,
\lambda)$ by $(\mu_{\lambda,c}, \mu_1,\theta).$
From Part 2 we have $k_{\mu_{\lambda,c}}(\theta)
=k_{\mu_{\lambda}}(\theta+ck_{\mu_{\lambda,c}}(\theta)).$
Denote $\alpha(\theta)=\theta+ck_{\mu_{\lambda,c}}(\theta)$ and therefore

$$\left[\begin{array}{c}\theta\\-k_{\mu_{\lambda,c}}(\theta)\end{array}\right]=\left[\begin{array}{cc}I_n&c\\0&\lambda\end{array}\right]\left[\begin{array}{c}\alpha(\theta)\\-k_{\mu_1}(\alpha(\theta))\end{array}\right]
$$
Here the function $\theta(\lambda)$ appearing in the statement of Proposition 3.3 is replaced by $\alpha(\theta).$
Infinite divisibility mentioned in Part 4 is obvious from Part 2. Also
$$\left[\begin{array}{c}\alpha\\-k_{\mu_{1}}(\alpha)\end{array}\right]=\left[\begin{array}{cc}I_n&-c\\0&\lambda\end{array}\right]\left[\begin{array}{c}\theta\\-k_{\mu_{\lambda,c}}(\theta)\end{array}\right]
$$ As a consequence $k_{\mu_{1}}(\alpha)
=k_{\mu_{1,c}}(\alpha-ck_{\mu_{1}}(\alpha)).$ If $\mu_1$ is a probability then $k_{\mu_{1}}(0)=0$ and clearly
$k_{\mu_{1,c}}(0)=0.$

\vspace{4mm}\noindent\textbf{Example: action of $H$ on the Poisson family.}
We take $n=1$ and $\mu_{\lambda}(dx)=\sum_{k=0}^{\infty}\frac{\lambda^k}{k!}$ meaning
$f(\lambda,k)=\frac{\lambda^k}{k!}1_{k\geq 0}.$ We get for $c>0$, $k\in \N$, and for the fixed parameter $x>0$
$$\mu_k(\lambda,c)=p_{k}(\lambda,c)1_{(0,\infty)}(cx)=e^{-\lambda-cx}\frac{(\lambda+cx)^{k-1}}{k!}.$$
The Laplace transform of $\mu_{\lambda}(dx)$ being $\exp(\lambda e^{\theta})$, implying $k_{\mu_{1}}(\theta)=e^{\theta}$, we get from Proposition 7.4
$k_{\mu(\lambda,c)}(\theta)=\lambda \exp(\theta +ck_{\mu(1,c)}(\theta))$. Writing $z(\theta)=k_{\mu(1,c)}(\theta))$ for simplicity we see that $z$ is obtained by the implicit equation
 $$z(\theta)=e^{\theta+z(\theta)}$$
or $z(\theta)=-W(e^{\theta}),$ where $W$ is a suitable branch of the Lambert function.
Actually the exponential family $F(\mu(\lambda,c))$ is rather familiar as the Abel family described in Letac and Mora (1990) formulas 4.15 and 4.36, obtained there by different road. To check this point we use the fact that the Poisson family has variance function $V(m)=m.$ Applying formula \eqref{TG} to $h_g(m)=m/(cm+\lambda)$ we get
$$V_{F(\mu(\lambda,c))}=\left(1+\frac{c}{\lambda}m\right)^2m.$$

\vspace{4mm}\noindent\textbf{Example: action of $H$ on the negative binomial family in $\R^n$.}

If $p=(p_1,\ldots,p_n)$ are positive numbers consider the Bernoulli measure in $\R^n$ defined by
$$\nu(dx)=\sum_{j=1}^np_n\delta_{e_j}(dx)$$ where $e_1,\ldots,e_n$ is the canonical basis of $\R^n$ as usual. The corresponding
negative binomial measure $\mu_{\lambda} $ is

$$\mu_{\lambda}(dx)=\sum_{i=0}^{\infty}\frac{(\lambda)_i}{n!}\nu^{*i}(dx)$$ where $(\lambda)_i=\lambda (\lambda+1)\ldots(\lambda+i-1)$  (and $(\lambda)_0=1$) is the Pochhamer symbol. Also $\nu^{*i}(dx)$ is the convolution power of $\nu^{*i}(dx)$ (and $\nu^{*0}(dx)=\delta_0(dx)$).
Because $\nu^{*i}(dx)$ is concentrated on the set $\{k\in \N^n; |k|=i\}$
the measure
$\mu_{\lambda}$ is well defined and its Laplace transform is
$$L_{\mu_{\lambda}}(\theta)=\frac{1}{(1-\langle p,e^{\theta}\>)^{\lambda}}$$ with the notation $\langle p,e^{\theta}\>=p_1e^{
\theta_1}+\cdots+p_ne^{\theta_n}.$ Here $\Theta(\mu_{\lambda})=
\{\theta\in \R^n; \langle p,e^{\theta}\><1\>\}$ and
$k_{\mu_{\lambda}}(\theta)=-\lambda \log (1- \langle p,e^{\theta}\>)$.
Let $c\in \R^n$ then from Proposition 7.4 we get
$$k_{\mu_{\lambda,c}}(\theta)=\lambda k_{\mu_{1}}(\theta+ck_{\mu_{1,c}}(\theta)).$$ For simplicity denote $z(\theta)=k_{\mu_{1,c}}(\theta).$ We get the implicit equation for the Laplace transform $L_{\mu_{1},c}(\theta)=e^{z(\theta)}:$

$$e^{-z(\theta)}=1-p_1e^{\theta_1+c_1z(\theta)}-\cdots-p_ne^{\theta_n+c_nz(\theta)}$$

\section{Conclusion} In this Part one,  we have gathered tools and results for the classification of the simple cubic exponential  families on $\R^n$. Their set  is denoted by $\mathcal{SP}_3(n)$.  This description  will be achieved in Part two. The main points are \begin{itemize}
\item The introduction of the group $G=G(n+1,\R)$ acting on the set of exponential families by $T_g$ introduced in \eqref{TG}, and its sub group $G_0$ of affine and Jorgensen tranformations;
\item  The  definition of the set of simple cubic exponential  families and its splitting into  $n+3$ $G$-orbits by Theorem 5.4;
\item The Lagrange formula in $\R^n$ and Theorem 6.2, which gives a necessary and sufficient condition for having an exponential  family concentrated on $\N^n$ from the consideration of its variance function.
\end{itemize}
Part two will perform the considerable work of describing the elements of each of the $n+3$ $G$-orbits, obtaining old and new distributions in $\R^n,$ mainly by using stopping times as it was done in the one dimensional case.

\end{document}